\documentclass[reqno]{amsart}
\usepackage{amscd}

\newcommand{\dv}{\operatorname{div}}

\newcommand{\A}{{\mathcal A}}

\renewcommand{\L}{{\mathcal L}}
\newcommand{\tm}{{\widetilde{M}}}
\newcommand{\G}{{\Gamma}}
\newcommand{\Gstar}{{{\stackrel{\ast}{G}}}}
\newcommand{\mustar}{{{\stackrel{\ast}{\mu}}}}

\newtheorem{thm}{Theorem}[section]
\newtheorem{lm}[thm]{Lemma}
\newtheorem{prop}[thm]{Proposition}
\newtheorem{crl}[thm]{Corollary}

\theoremstyle{definition}
\newtheorem{df}[thm]{Definition}

\theoremstyle{remark}
\newtheorem{rem}[thm]{Remark}

\begin{document}

\title{Martin points  on open manifolds of non-positive curvature}
\author{Jianguo Cao}
\address{Department of Mathematics, University of Notre Dame, Notre Dame, IN 46556, USA} \email{jcao@nd.edu}
\author{Huijun Fan}
\thanks{The first author is supported in part by an NSF grant. The second author was partially
supported by the Research Fund for returned overseas Chinese
Scholars 20010107 and by the NSFC for Young Scholars 10401001}
\address{School of Mathematical Sciences, Peking University, Beijing 100875, China
and Max-Planck Institute for Mathematics, Inselstr. 22-26, 04103,
Leipzig, Germany } \email{fanhj@math.pku.edu.cn}
\author{Fran\c{c}ois Ledrappier}
\address{Department of Mathematics, University of Notre Dame, Notre Dame, IN 46556, USA} \email{fledrapp@nd.edu}

\maketitle

\begin{abstract}

The Martin boundary of a Cartan-Hadamard manifold describes a fine geometric structure at infinity, which is a sub-space of positive harmonic
functions. We describe conditions which ensure that some points of the sphere at infinity belong to the Martin boundary as well. In the case of
the universal cover of a compact manifold with Ballmann rank one, we show that Martin points are generic and of full harmonic measure. The
result of this paper provides a partial answer to an open problem of S. T. Yau.

\end{abstract}

\numberwithin{equation}{section}
\section{Introduction}
Let $\tm$ be a Cartan-Hadamard manifold, a simply connected Riemannian manifold with nonpositive curvature. Then,  $\tm $ is  homeomorphic to an
open ball, and there are two natural compactifications of $\tm $ associated to the metric.

 Fix $x_0 \in \tm$. For $z \in \tm $, define the continuous function $b_z $ on $\tm $ by:
$$ b_z (x) \; = \; d(x,z) - d(x_0, z), $$
where $d$ denotes the Riemannian distance on $\tm$. The functions
$b_z, z \in \tm $ are equicontinuous and uniformly bounded on
compact subsets of $\tm$. They form a relatively compact set of
functions for the topology of  uniform convergence on compact
sets. The closure of $\{ z \mapsto b_z \}$ is the geometric
compactification of $\tm$. Let $\tm (\infty)$ be the boundary of
$\tm $ in its geometric compactification. The set $\tm (\infty) $,
endowed with the relative topology, is homeomorphic to a sphere.
Let $T(\tm )$ be the tangent bundle of $\tm$,  and  $ S_{x}\tm =
\{\vec v \in T_x(\tm) | \quad \| \vec v \| = 1  \}$ be the unit
tangent sphere of $\tm $ at $x$.
 For any $x \in \tm$, the map
$P_{x} : S_{x}\tm  \mapsto \tm (\infty)$  which associates to $v \in S_{x}\tm $ the  point $P_{x} v = \sigma _v (+\infty)$ realizes this
homeomorphism, where $\sigma _v$ is the geodesic with initial condition  $v$ and, for a geodesic $\sigma$  in $\tm$, we denote
$\sigma (\pm \infty ) = \lim _{t \to +\infty } \sigma (\pm t)$ the
corresponding  points  of $\tm (\infty) $.

Assume that $\tm $ admits a Green function  $G(.,.)$ for  the Laplace operator. For $z \in \tm $, define the continuous function $h_z $ on $\tm
$ by:
$$ h_z (x) \; = \; \log G(x,z) - \log G(x_0, z). $$
By Harnack inequality, the functions $h_z, z \in \tm $ are equicontinuous and uniformly bounded on compact
subsets of $\tm$ not containing $z$. They form a relatively compact set of functions for the topology of  uniform convergence on compact sets.
The closure of $\{ z \mapsto h_z \}$ is the Martin compactification of $\tm$.

 For Euclidean spaces, the Martin
   compactification  is reduced to the Alexandroff one-point compactification.  If the  sectional curvatures of $\tm$ are pinched between two negative constant, then the Martin compactification coincide with the
    geometric compactification \cite {AS}. In general, the presence of flats amidst  negative curvature is a source of more intricate Martin compactification:   for symmetric spaces, the Martin compactification
   has been described in \cite {GJT} and is a nontrivial continuous extension of the geometric compactification;  the general description of the Martin compactification  of a product
  is not known in general,  see \cite {MV} and the references therein for the latest results. In these two cases, every geodesic belongs to a flat space. It is believed that, if many geodesics are not within a totally geodesic flat subspace,  then the Martin compactification is the geometric  compactification. See \cite {Ba4} for a first example.

Following  Ancona's programme (see \cite{An}), the same discussion applies to  the general uniform elliptic operator $\L$ of second order
in a general Cartan-Hadamard  manifold $\tm$ of dimension $n\ge 2$ with bounded geometry. The elliptic operator $\L$ has the
form
$$
\L(u):=\dv(\A(\nabla u))+B\cdot\nabla u+\dv(uC)+\gamma u.
$$
The conditions for the coefficients will be given in the next
section. If $\L u=0$, then $u$ is called a $\L$-harmonic function.
We still denote by $G(\cdot,\cdot)$ the Green function of $\L$ and
define $h_z(x)$ as before. Again fix $x_0 \in \tm$.

\begin{df}[Poisson kernel function]\label{Poisson} A Poisson kernel function
$k_\xi(x)$ of $\L$ at $\xi\in \tm(\infty)$ is a positive
$\L$-harmonic function on $\tm$ such that:
\begin{equation} k_\xi (x_0 ) = 1, \; k_\xi (y)=O(G_{x_0}(y)) \; \text { as }\; y \to \xi' \not = \xi, \end{equation}
\end{df}
\begin{df}[Martin point]
 We say that a point $\xi \in \tm(\infty) $ is a Martin point of $\L$ if it satisfies the following properties:
\begin{itemize}
\item a) There exists a Poisson kernel function  $k_\xi $ of $\L$
at $\xi$, \item b) the Poisson kernel function  is unique, and
\item c) if $y_n \to \xi$, then $h_{y_n} \to \log k_\xi$ uniformly
on compact sets.
\end{itemize}
\end {df}

In this paper, we want to describe  Martin points of $\L$ for
Cartan-Hadamard manifolds. For that purpose, we introduce several local notions of negative curvature along a geodesic in $\tm$.  For a vector $v \in S\tm$, the {\it rank } of $v$ is the dimension of the space of parallel Jacobi fields along the geodesic $\sigma _v $ with initial condition $v$. Clearly, $1 \le \text { rank } v \le \text { dim } \tm$. The geodesic  rank of the manifold $\tm$ is the minimum value of $\{ \text {rank }v, v \in S\tm \}$. For locally symmetric spaces, the geodesic rank coincide with the real rank of the real  algebraic group of isometries of $\tm$.

A geodesic $\sigma $ is called {\it rank one}  if rank of $\sigma '(0) $ is equal to $1$. A geodesic in $\tm$ is called {\it regular }  if it
does not bound a totally geodesic flat half-space. Rank one  geodesics are
regular. In the next section, we introduce the notion of { \it hyperbolic} geodesic in $\tm$.
It is a precise qualitative property which expresses that the geodesic has an infinite number of segments surrounded
by enough negative curvature. Geodesics in flats, or even geodesics converging to flats are not hyperbolic.
Our main result is:

\begin {thm}\label{theo-axis} Let $\tm$ be a Cartan-Hadamard
manifold with bounded geometry, $\L$ a uniformly elliptic, weakly coercive and bounded  second order operator and $\sigma:{\mathbb R}\mapsto \tm$  a hyperbolic geodesic. Then $\sigma(+\infty)$ is a Martin point
of $\L$. In particular, if the Laplace operator $\Delta $  is weakly coercive,  $\sigma(+\infty)$ is a Martin point for  $\Delta $.
\end {thm}

An {\it axis} in $\tm$ is a geodesic which is invariant  by an   isometry of $\tm$ with two fixed points at infinity.  We will see that regular axes are hyperbolic.
\begin{crl}\label{crl-axis}  Let $\tm$ be a Cartan-Hadamard
manifold with bounded geometry, $\L$ a uniformly elliptic, weakly coercive and bounded  second order operator and $\sigma:{\mathbb R}\mapsto \tm$  an axis  such that ${\sigma}$ is not  a boundary of
any totally geodesic half-plane. Then ${\sigma}$ is hyperbolic and $\sigma(+\infty)$ is a Martin point
of $\L$.
\end{crl}

\begin{rem} If the sectional curvature of $\tm$ is pinched, then
Ancona(\cite{An}) has proved that the Martin boundary
$\partial_\L\tm$ of $\tm$ with respect to $\L$ is homeomorphic to
the geometrical boundary $\tm(\infty)$. Our result extends Ancona's results  to nonpinched manifold, at least at extremities of hyperbolic geodesics.
\end{rem}

\

 In the rest of the paper, we show that if $\tm$ is rank one and admits a cocompact group of isometries, then there are many hyperbolic geodesics. So assume that the manifold $\tm $ is  the universal cover of a compact manifold $M$. Then, $\tm$ has bounded geometry as soon as the metric is of class $C^3$, and the Laplace operator admits a Green function as soon as $M$ is not a  2-dimensional torus. Moreover, the geodesic rank rigidity results of Ballmann \cite {Ba2} and Burns-Spatzier \cite {BS} asserts that $\tm $ can
 be written uniquely as a product of Euclidean spaces, symmetric spaces and the universal covers of  rank one spaces
 (see \cite {K1}, Appendix,  for the existence of a cocompact action on the third  factors). We shall therefore concentrate on rank one manifolds.  We have:

\begin{crl}\label{main-crl} Let $\tm$ be the universal cover of a compact
Riemannian manifold of class $C^3$,  non-positive curvature and geodesic  rank 1,   $\L$ a uniformly elliptic, weakly  coercive and bounded  second order operator on $\tm$ and $\sigma:{\mathbb R}\mapsto \tm$  a regular axis. Then,  $\sigma(+\infty)$ is a Martin point
of $\L$. In particular, Martin points are dense in $\tm (\infty)$.
\end{crl}

\

Let $\G = \pi_1(M)$ be the covering group. Recall that the action of $\G$ by isometries on $\tm$ extends  to a continuous action on $\tm (\infty)$.
We set  $(\tm (\infty ) \times \tm (\infty ))^*$ for the set of pairs of distinct points in $\tm (\infty)$.
We say that a finite positive measure $\mu$ on $\tm (\infty)$ is {\it geodesic ergodic } if
\begin {itemize}
\item 1) The support of the measure  $\mu\times \mu $ is $(\tm (\infty ) \times \tm (\infty ))^*$.
\item 2) For $\mu \times \mu $ almost every $(\eta, \xi)$, there is a unique  geodesic $\sigma _{\eta,\xi}$ such that: $\sigma _{\eta,\xi}(-\infty ) = \eta, \sigma _{\eta,\xi} (+\infty ) = \xi$, and  $\sigma _{\eta,\xi}$ is rank one.
\item 3) The measure $\mu \times \mu $ is $\G$ quasi-invariant and ergodic: the diagonal action of   $\G$ preserves the $(\mu \times \mu )$-negligible
subsets  of  $(\tm (\infty ) \times \tm (\infty ))^*$;  and all $ \G$-invariant measurable subsets of $(\tm (\infty ) \times \tm (\infty ))^*$
 are either negligible or co-negligible.
\end {itemize}

Examples of geodesic ergodic measures are the Patterson-Sullivan measure (see \cite {K1}), other Gibbs measures constructed along the same lines, the harmonic measure for the Laplace operator on $\tm$ (see \cite {BL}), or analogously other harmonic measures associated to  Markov equivariant symmetric operators on $\tm$ or on $\G$ (\cite {Ka}). It is not known, even for surfaces,  whether the visibility measure, obtained by projecting under $P_{x_0}$  the Lebesgue  measure of the  sphere $S_{x_0}\tm$ is geodesic  ergodic. We have

\begin{thm}\label{many-mart} Let $\tm$ be the universal cover of a compact
Riemannian manifold of class $C^3$, non-positive curvature and geodesic  rank one, and   $\L$ a uniformly elliptic, weakly  coercive
and bounded second order operator on $\tm$. Then the set of Martin points is a generic subset of $\tm (\infty)$: it
contains a countable intersection of open dense subsets. Moreover, the set of Martin points has full measure for any
geodesic ergodic measure.
\end{thm}

In the next section, we introduce the necessary definitions and present the general scheme of the proofs.
In section 3, we recall the potential theory of weakly coercive operators, and  section 4 contains the geometric properties
of  hyperbolic geodesics  we shall use. Theorem \ref{theo-axis} reduces to Propositions we prove in section 5,
and Theorem \ref{many-mart} is proven in  section 6.

\

\section{Precise statements of results and strategy of the proofs.}\

Let $\tm$ be a complete Riemannian manifold of dimension $n\ge 2$.
If $d(y, z)$ is sufficiently small, we let $\mathbb P_y^z$ denote
the  parallel transport from $y$ to $z$ along the unique
length-minimizing geodesic segment.

We say  that $\tm$ has {\it bounded geometry}  if
there exists $r_0>0$ such that for any ball $B(x,r_0)\subset \tm$,
there exists a chart $\chi:B(x,r_0)\rightarrow {\mathbb R}^n$
satisfying a uniform {\it first order}  quasi-isometry condition:
\begin{align}\label{boud-geom}
& C_0^{-1}d(y,z)\le ||D\chi_x|_y -D\chi_x|_z ||^* \le C_0 d(y,z),\nonumber\\
&\forall y,z\in B(x,r_0), \text{with a constant $C_0$ independent
of $x$},
\end{align}
where
$$
 ||D\chi_x|_y-D\chi_x|_z||^* = ||\chi_x(y)-\chi_x(z)|| + \max_{ \|\vec v \| = 1 }\{ \|(D\chi_x) |_y \vec v -
 (D\chi_x) |_z (\mathbb P_y^z \vec v ) \|
 \}. $$

If   the derivative of the
curvature of $\tm$ is bounded, and if  the injectivity radius of $\tm>0$,
then (\ref{boud-geom}) holds.\

Consider the following elliptic operator $\L$:
\begin{equation}
\L(u):=\dv(\A(\nabla u))+B\cdot\nabla u+\dv(uC)+\gamma u,
\end{equation}
where $\A$ is a section of $End(T\tm)$, $B$ and $C$ are vector
fields on $\tm$ and $\gamma$ is a function.\

\begin{df} The operator $\L$ is called uniformly elliptic if there is $\lambda >1$ such that:
\begin{equation}\label{elli-con1}
\forall (x,u)\in T\tm, \lambda^{-1}||u||^2\le
\langle\A_x(u),u\rangle\le \lambda ||u||^2.
\end{equation}
\end{df}

\begin{df}  The operator $\L$ is said to be bounded  if  there is $\lambda >0$ such that:
\begin{equation}\label{elli-con2}
\forall x\in \tm, ||B||_{L^\infty(B(x,r_0))},
||C||_{L^\infty(B(x,r_0))}, ||\gamma||_{L^\infty(B(x,r_0))}\le
\lambda.
\end{equation}
\end{df}

\begin{df} The function $G:\tm\times \tm\mapsto (0,+\infty)$ is called a Green
function of $\L$, if $G$ is continuous, and for any $x\in \tm$,
$G_x(y):=G(y,x)$ is a $\L$-potential on $\tm$ and is $\L$-harmonic
on $\tm\setminus \{x\}$ such that
$$
\L(G_x)=-\delta_x.
$$
\end{df}

\begin{df} The operator $\L$ is called weakly coercive, if there exists
$\epsilon>0$ and a positive superharmonic function on $\tm$ with
respect to the operator $\L+\epsilon I$.
\end{df}
So if $\L$ is weakly coercive for some $\epsilon>0$, then for any
$0\le t<\epsilon$, the operator $\L+tI$ has a Green function
$G^t$. \

\

Let now $\tm $ be a Cartan-Hadamard manifold, and for  $\sigma:{\mathbb R}\mapsto \tm$  a geodesic line of unit
speed, set
$$
\mathcal{U}_h(\sigma({\mathbb R}))=\{y|d(y,\sigma({\mathbb R}))=h\}.
$$
Since $\sigma(\mathbb R)$ is a closed convex subset of $\tm$, there is the nearest-point projection: $\mathcal{P}_\sigma: \tm \mapsto
\sigma(\mathbb R)$.
We define
$$
S^\perp_h(\sigma(t))=\mathcal{P}_\sigma^{-1}(\sigma(t)) \cap
 \mathcal{U}_h(\sigma({\mathbb R})),
$$
and
$$
\eta_{\sigma([t_1,t_2])}(h)=d_{\mathcal{U}_h}(S^\perp_h(\sigma(t_1)),
S^\perp_h(\sigma(t_2))),
$$
where $d_{\mathcal{U}_h}(\cdot,\cdot)$ is the distance
function of the Riemannian hypersurface $(\mathcal{U}_h,
g|_{\mathcal{U}_h})$.

\

The following notion is a  way of expressing at a finite distance that the geodesic $\sigma $ does not bound a flat half space:
\begin{df} A geodesic  $\sigma:{\mathbb
R}\mapsto \tm $ is said to be $(h,T,\delta)$-non flat at $t$ if we have:
$$ \eta_{\sigma([t,t+T])}(h) \; > \; T + \delta h.$$
\end{df}

Properties of $(h,T,\delta)$-non flat geodesics are recalled in Section 4. In particular, by Proposition \ref {nonflat-cont}
there exists a number $\varepsilon ^* = \varepsilon ^* (\tm, h ,T)$ such that if the geodesic $\sigma $ is $(h,T, {\pi}/{2})$-non flat at $0$, and $\tau $ is another geodesic satisfying
$$ \tau (0) \;= \; \sigma (0) \; \text{and} \; \angle_{\sigma(0)}(\tau'(0),\sigma'(0)) \; <\; \varepsilon ^*,$$
then the geodesic $\tau $ is $(h,T,  {\pi}/{4})$-non flat at $0$.

Let us now choose $\varepsilon ^* <  {\pi}/{4}$ and set
$$ T_1 \; = \; T_1 (\tm,h,T)\; = \; T + \frac {h} {\tan \varepsilon ^*}.$$

\begin{df}\label{bar} We say that the geodesic $\sigma $ admits a $(h,T,R) $ barrier
 if there exist $t_i, i= 1, 2,\dots, 6$ with $T_1 < t_{i+1} -t_i< T_1 + R$
 and $t_3 + T <0 < t_4$ such that the geodesic $\sigma $ is $(h,T, {\pi}/{2})$-non flat at
 $t_i$, for $i = 1,2,\dots, 6$.
\end{df}

\begin{rem} \label{reverse} Observe that if a  geodesic $\sigma $ is $(h,T, \pi/2)$-non
flat at 0, the geodesic $-\sigma $ obtained by reversing time is  $(h,T, \pi/2)$-non flat
at $-T$. Consequently, if the geodesic $\sigma $ admits a $(h,T,R)$ barrier, the geodesic
 $-\sigma $ admits a $(h,T,R)$ barrier as well, with $t'_i = -t_{7-i} -T$.
\end{rem}

\begin{df}\label{hyp} We say that the geodesic $\sigma $ is hyperbolic if there are $h,T,R $ and a sequence $t^*_i \to + \infty$ such that $\sigma (\cdot -t^*_i)$ admits a $(h,T,R)$ barrier.
\end{df}

We have defined all elements of Theorem \ref{theo-axis} that we recall:

{\bf Theorem \ref{theo-axis}} {\it  Let $\tm$ be a Cartan-Hadamard
manifold with bounded geometry, $\L$ a uniformly elliptic, weakly coercive and bounded
second order operator and $\sigma:{\mathbb R}\mapsto \tm$  a hyperbolic geodesic.
Then $\sigma(+\infty)$ is a Martin point
of $\L$.}

\

In order to prove Theorem \ref {theo-axis}, we define the
families of cones
$$
\Gamma_{\sigma,t,\theta}=\{x\in\widetilde{M}|\angle_{\sigma(t)}(\sigma'(t),x)<\theta\}.
$$

\begin{thm}\label {barrier}
Suppose the geodesic $\tau $  admits a $(h,T,R) $ barrier. Set $T_2 = 3(T_1 +R) +T$. Then there is a constant $C = C(\tm, h,T,R)$ such that the Green function $G(x,y)$ satisfies
 \begin{align}\label{unif-G}  &G(x,y) \le C G(x, \tau(0))G(\tau (0),y),\\
&\forall x \in \tm \setminus \Gamma _{ {\tau},-T_2, \pi/2},\;  \forall y \in \Gamma _{{\tau}, T_2, \pi/2}\nonumber
\end{align}
 and the Green function $g(x,y)$ in $\tm \setminus \Gamma_{ {\tau},2T_2,\pi/2} $
satisfies
\begin{align}\label{unif-cons}
&g(x,y)\le C g(x,\tau(0))
g(\tau(0),y),\\
&\forall x \in \tm \setminus\in\Gamma_{
{\tau},-T_2, \pi/2},\; \forall y \in  \Gamma_{ {\tau},T_2 ,\pi/2} \setminus \Gamma_{
{\tau}, 2T_2,\pi/2}.\nonumber
\end{align}
\end{thm}

 Recall Definition \ref {Poisson} of a Poisson kernel function, and call $C_\xi$ the cone of functions
 positively proportional to a Poisson kernel function  at $\xi \in
 \tm(\infty)$. Then,
\begin{prop}\label{unicity} Assume $\tau $ is a hyperbolic geodesic with $\xi = \tau (+\infty) $. Then,  $\text {dim}\;
C_{\xi}\le 1$.
\end{prop}

\begin{prop}\label {existence}  Assume $\tau $ is a hyperbolic geodesic with $\xi = \tau (+\infty) $, and consider the functions $ k_z (x)  =  \frac {G(x,z)}{ G(x_0, z)}. $ Then, if $k_\xi$ is a limit point of $k_z$ as $z \to \xi$, $k_\xi \in C_{\xi}.$
\end{prop}

Theorem \ref{theo-axis} follows directly from Propositions \ref{unicity}  and \ref {existence}. In Section 5, we prove Theorem \ref {barrier} and explain how Propositions \ref{unicity}  and \ref {existence} follow from Theorem \ref{barrier}. In \cite {An}, (\ref {unif-G}) is called the Boundary Harnack Inequality and is a key step in the proof. For establishing (\ref{unif-G}), our task is to use as little negative curvature as we find it necessary. The proof follows the ideas from \cite{An}, but given the delicate arguments involved, we prefer writing it in whole detail. Then, following \cite{An}'s scheme, Propositions \ref{unicity}  and \ref {existence} follow from Theorem \ref{barrier}. Our observation is that it is sufficient to have an infinite number of disjoint barriers converging to $\xi$, not necessarily a uniform estimate everywhere. Again we write the detailed  proof for the sake of completeness.

\

Assume now that  $\tilde{\sigma}:{\mathbb R}\mapsto \widetilde{M}$ is
an axis and suppose that $\tilde{\sigma}$ is not the
boundary of any totally geodesic half plane. Then, there exist $h_0$ and $\delta_0$ such that for any $k\in {\mathbb N}$,
 there is an integer $n$ such that $\tilde \sigma $ is $(h_0,n L, k\delta_0)$-non flat at $0$, where  $L$
is the period of axis $\tilde{\sigma}$.

 Indeed, since $\tilde{\sigma}$ is invariant by an isometry,
$\tilde{\sigma}$ is not the boundary of any totally geodesic flat
 two-dimensional quarter. Thus, by corollary \ref{flat} there exist $T_0, h_0$ and $\delta_0$ such that
$$
\eta_{\tilde \sigma([0,T_0])}(h_0)-T_0\ge \delta_0>0.
$$
Choose $n_0> T_0/L $ to be an integer. Thus, since the function $T
\mapsto \eta_{\tau [0,T]}(h_0) - T $ is nondecreasing (see
Proposition \ref {subad} and Proposition \ref{prop3.23}(5)):
$$
\eta_{\tilde \sigma([0,n_0L])}(h_0) -n_0L \ge \delta_0.
$$

For any integer $k$, we get, using semiaddivity (\ref{subad}) and the periodicity of $\tilde \sigma$
$$
\eta_{\tilde{\sigma}([0,kn_0L])}(h_0) -kn_0L \ge k \delta_0,
$$
which is the desired property by setting $n = n_0 k$.

By invariance under isometries, the axis  $\tilde \sigma $ is also
$(h_0,n L, k\delta_0)$-non flat at $KL $, for all $K \in \mathbb
N$.   By choosing $k$ such that $k\delta _0 > \frac {\pi}{2}h_0$,
and $t_i, i= 1,2,\dots, 6$ also multiples of $L$, we find a number
$R$ such that the axis $\tilde \sigma $ admits a $(h_0, nL, R)$
barrier. By invariance by isometries again, the axis  $\tilde
\sigma$ is a hyperbolic geodesic. Corollary \ref{crl-axis} is
therefore a particular case of Theorem \ref{theo-axis}.

\

Consider   the case when the Cartan-Hadamard manifold $\tm$ is the universal cover of a compact Riemannian  manifold $M$  of geodesic rank one. Set $SM$ for  the unit tangent bundle of $M$. A unit tangent vector  $v \in SM$ is said to be regular, $(h,T,\delta)$-non flat, admitting a $(h,T,R)$ barrier or hyperbolic if any geodesic $\sigma _{\tilde v}$ defined by a lift $\tilde v$ of $v$ to $S\tm$ has the same property.  Ballmann (\cite {Ba1}) showed that unit tangent vectors to regular closed geodesics are dense in $SM$. Therefore Corollary  \ref{main-crl} directly follows from Theorem \ref{crl-axis}. The geodesic flow is a one parameter group $\varphi_t, t\in {\mathbb R}$ of diffeomorphisms of $SM$. There is a unique $\varphi $-invariant probability measure $\bar \nu$  on $SM $  which realizes the topological entropy. The measure $\bar \nu$ has full support on $SM$  and the geodesic flow is ergodic for $\bar \nu$ (\cite {K2}). Therefore:

\begin{prop}\label{dense} Let $M$ be a compact Riemannian manifold of nonpositive sectional curvature and geodesic rank 1. Then the set of hyperbolic unit tangent vectors contains a countable  intersection of open dense sets in $SM$. Moreover, it has full measure for $\bar \nu$.
\end{prop}
\begin{proof} We know that  a unit tangent vector to a regular closed  geodesic admits a $(h,T,R)$ barrier for some $h,T$ and $R$. By proposition \ref{nonflat-cont}, there is an open neighborhood $\mathcal O$ of such a unit vector $v$ such that all  $v' \in {\mathcal O}$ also admit a $(h,T,R)$ barrier. Since the measure $\bar \nu$ is ergodic and has full support, for all positive $K$ the set ${\mathcal O}_K$ of $v \in SM$ such that the geodesic ray $\sigma _v([K,\infty))$ intersects $\mathcal O$ is open dense in $SM$ and has full $\bar \nu$ measure. The set $\cap_K {\mathcal O}_K$ is a countable intersection of open dense sets of full $\bar \nu$ measure. By definition, any unit vector in $\cap_K {\mathcal O}_K$ is hyperbolic.
\end{proof}

To prove Theorem \ref{many-mart}, we still  have to verify that the large set of unit vectors of Proposition \ref{dense}  lifts and projects to a large subset of $\tm (\infty)$. This relies on the properties of the measure $\bar \nu$ which have been established in \cite {K2}, see Section 6.

\begin{rem} In the case when $\tm$ is the universal cover of a compact rank 1 manifold, the Laplace operator $\Delta $ is weakly coercive (see below section 3) and clearly uniformly elliptic and bounded. The conclusions of Corollary  \ref{main-crl} and Theorem \ref{many-mart} hold for $\L = \Delta$.
\end {rem}

\

\section{Preliminaries (Elliptic operators, Green functions and their estimates)}\

Let $\tm$ be a Cartan-Hadamard manifold with bounded geometry and
$\L$ a uniformly elliptic, weakly coercive and bounded  second
order operator. Let $\mu$ be a positive measure on $\tm$. Define
$G\mu(x):=\int_{\tm} G(x,y)d\mu(y)$. If $G\mu$ is not identically
$+\infty$, $G\mu$ is the only potential satisfying
$\L(G\mu)=-\mu$.\

There are two important estimates (see \cite{An}):\
\begin{itemize}
\item For each $\omega=B(x,r_0)\subset \tm $ , and every $t, 0\le
t\le 1$, the Green function $g^t$ related to $\L+tI$ over $\omega$
satisfies
\begin{equation}
g^t(y,z)\ge C, \forall y, z\in B(x,r_0/2), \;\text{and}\;
g^t(y,z)\le C^{-1}, \;\text{if}\;d(y,z)\ge \frac{r_0}{4},
\end{equation}
where $C=C(\L)$ is independent of $x$ and $t$. \

\item (Harnack inequality) If $u>0$ is a $\L+tI$-harmonic function
on $B(x,r_0)$, then
\begin{equation}
C^{-1}u(x)\le u(y)\le Cu(x),
\end{equation}
where $C=C(\L)>0$.
\end{itemize}
The adjoint operator $\L^*$ of $\L$ is given by the formula:
$$
\L^*(u)=\dv(\A^*(\nabla u))-\dv(B\cdot u)-C\cdot \nabla u+\gamma
u.
$$

Note that the Green function $\Gstar(x,y)$ of $\L^*$ satisfies
$\Gstar(x,y)=G(y,x)$. \

\begin{lm}[\cite{An}, Lemma 1]\label{anco-lm1} For each positive measure $\mu$ on $\tm$ and each $t,
0\le  t<\epsilon$, we have
$$
G^t(\mu)=G(\mu)+G(G^t(\mu)).
$$
\end{lm}

Let $\mustar_x$ be the $\L^*$-harmonic measure of a point $x\in
\Omega$,  where $\Omega$ is a bounded region in $\tm$. We have

\begin{lm}[\cite{An}, Lemma 3] Let $g$ be the $\L$-Green function of
$\Omega$, and let $g_x(y)=0$ for $y\not\in\Omega$, then
$$
\L(g_x)=-\delta_x+\mustar_x.
$$
\end{lm}

\begin{proof} We have the representation formula of $g(x,y)$ in
terms of $G(x,y)$ and the harmonic measure $\mustar_x$:
$$
g_x(y)=\Gstar_y(x)-\int_{\partial \Omega}\Gstar_y(z)
d\mustar_x(z),
$$
for $x\in \Omega, y\in \tm$. Then we have
$$
g_x=G_x-G(\mustar_x),
$$
and so
$$
\L(g_x)=-\delta_x+\mustar_x.
$$
\end{proof}

Denote by $g^t$ the $\L+tI$-Green function of $\Omega$, and by
$\stackrel{\ast}{\mu}{\!}^t_x$ the $\L^*+tI$-harmonic measure of
$x$ in $\Omega$, we have

\begin{lm}[\cite{An}, Lemma 4]\label{Gree-meas} If $0\le t<\epsilon, x\in \Omega$ and
$g_x\le k g^t_x$ for some $k>0$ and outside some compact subset of
$\Omega$, then we have
$$
\mustar_x\le k \stackrel{\ast}{\mu}{\!}^t_x.
$$
\end{lm}

\begin{df} Let $\Omega$ be a not necessarily bounded region in
$\tm$. Let $x\in \Omega$, the "reduit" of $G_x$ on
$\overline{\Omega}^c$ is defined as
$$
R^{\overline{\Omega}^c}_{G_x}:=\inf\{s|s>0\;\text{is}\;
\L-\text{superharmonic on }\;\tm, \text{and}\;s\ge
G_x\;\text{on}\;\overline{\Omega}^c\}.
$$
This reduit is an $\L$-potential, and if we put
$\nu_x=-\L(R^{\overline{\Omega}^c}_{G_x})$, then $\forall z\in \tm
\setminus \overline{\Omega}$, we have the formula:
$$
\Gstar_z(x)=G_x(z)=\int G(z,y)d\nu_x(y),
$$
where $\nu_x$ is supported by $\partial\Omega$.
\end{df}

\begin{prop}[\cite{An}, Proposition 7]\label{anco-prop-7} There is a constant
$C=C(\tm,\lambda,\epsilon)>0$ such that if $x,y\in \tm$ and
$d(x,y)=1$, then
$$
\frac{1}{C}\le G^t(x,y)\le C,\;\text{for}\;0\le t<\epsilon
$$
\end{prop}

\begin{lm}[\cite{An}, Lemma 9]\label{deca-1} There exists a constant
$\delta=\delta(\tm,\lambda,\epsilon), 0<\delta<1$, such that for
each ball $B(x,1)$ in $\tm$, the $\L$-harmonic measure $\mu_x$ of
$x$ in $B(x,1)$ and the similar $\L+\epsilon I$ harmonic measure
$\mu^\epsilon_x$ satisfy
$$
\mu_x\le (1-\delta)\mu^\epsilon_x.
$$
\end{lm}

\begin{prop}[\cite{An}, Proposition 10]\label{decay-2} There are positive numbers $C$ and
$ \alpha$ such that
\begin{equation}
G(x,y)\le C e^{-\alpha d(x,y)}G^\epsilon(x,y), \forall x,y\in \tm,
\end{equation}
where $C$ and $\alpha$ depend only on $\tm,\lambda$ and $\epsilon$.
\end{prop}

\begin{proof} By induction on $k\in {\mathbb N}, k\ge 1$, we prove
that $G(x,y)\le (1-\delta)^{k-1}G^\epsilon(x,y)$, for $d(x,y)=k$
and $\delta$ which is given by Lemma \ref{deca-1}.\

When $k=1$, we have $G(x,y)\le G^\epsilon(x,y)$, since
$G^\epsilon(x,y)$ is a $\L$-superharmonic function.\

Assume that the inequality holds for $d(x,y)=k$. We want to prove
that it holds for $d(x,y)=k+1$. By maximum principle, one has
$$
G_x(z)\le (1-\delta)^{k-1}G^\epsilon_x(z), \forall z\in \tm\setminus
B(x,k).
$$
In particular, for $z\in \partial B(y,1)$. Hence
$$
G_x(y)=\int_{\partial B(y,1)}G_x(z)d\mu_y(z)\le
(1-\delta)^{k-1}\int G^\epsilon_x(z) d\mu_y(z).
$$
Now by Lemma \ref{deca-1},
$$
G_x(y)\le (1-\delta)^{k}\int G^\epsilon_x(z)
d\mu_y^\epsilon(z)=(1-\delta)^k G^\epsilon_x(y).
$$
This proves the proposition for $d(x,y)$ being integer. The
general case follows by the fact $G_x\le G^\epsilon_x$ and Harnack
inequality for $G^\epsilon_x$.
\end{proof}

\begin{rem}\label{deca-rema} Let $\Omega=B(x,r)$. Then proposition \ref{decay-2}
holds for $G_\Omega$ and $G^\epsilon_\Omega$, with the constants
$C,\alpha$ independent of $r$. This is because if we proved the
estimate for $d(x,y)\le r-1$. Then if $r-1\le d(x,y)<r$, by
maximum principle, we have
$$
G_\Omega(x,y)\le C e^{-\alpha(r-1)}G^\epsilon_\Omega(x,y).
$$
\end{rem}

\begin{rem}\label{deca-rema2} By Harnack inequality and Proposition \ref{anco-prop-7}, it is easy to
obtain the lower bound estimate of $G(x,y)$:
$$
ce^{-\beta d(x,y)}\le G(x,y),
$$
where $c, \beta>0$ only depend on the bounded geometry of $\tm$
and the operator $\L$.
\end{rem}

\begin{crl}[\cite{An}, Corollary 11]\label{comp-meas} Given $\delta>0$, there exists
$R=R(\tm,\lambda, \epsilon,\delta)$ such that $\forall x\in \tm$, and
$\forall r\ge R$, the $\L$-harmonic measure $\mu_x$ of $x$ in
$B(x,r)$ and the similar $\L+\epsilon I$ harmonic measure
$\mu^\epsilon_x$ satisfy:
\begin{equation}
\mu_x\le \delta\mu_x^\epsilon.
\end{equation}
\end{crl}

\begin{proof} For given $\delta>0$, we can find
$R=R(\tm,\lambda,\epsilon,\delta)$ such that $C e^{-\alpha
d(x,y)}\le \delta$ for $y$ near $\partial B(x,r)$ for any $r\ge
R$, where $C$ and $\alpha$ are from Remark \ref{deca-rema}. So
$$
G_B(x,y)\le \delta G_B^\epsilon(x,y), \;\text{for}\; y
\;\text{near}\; \partial B(x,r),
$$
i.e.,
$$
\Gstar_{B,x}\le \delta \stackrel{\ast}{G}{\!}^\epsilon_{B,x}.
$$
By Lemma \ref{Gree-meas}, we have
$$
\mu_x\le \delta\mu^\epsilon_x.
$$
\end{proof}

\

\vspace{5mm}

\

Assume now that  the Cartan-Hadamard manifold $\tm$ is cocompact, i.e.,
it is the universal cover of some compact Riemannian manifold $M$
with the lifted metric. Furthermore, we assume $M$ is of geodesic
rank 1. It is known that the fundamental group $\pi_1(M)$ of $M$
contains a free group $F_2$, and hence $\pi_1(M)$ is non-amenable.
By Brooks's result, the first eigenvalue of Laplace operator $$
\lambda_1(\tm)=\inf_{f\in H^{1,2}(\tm)}\frac{\int_\tm|\nabla
f|^2}{\int_\tm |f|^2}>0.
$$

Now let $G(x,y)$ be the Green function of the Laplace operator
$\Delta$ on $\tm$. Since $\tm$ is cocompact, the sectional
curvature $|K_{\tm}|$ and its derivative are  bounded and the injectivity radius
$inj(\tm)$ is positive. Thus $\tm$ has the ``bounded geometry" property
(\ref{boud-geom}). On the other hand, Laplace operator $\Delta$
satisfies (\ref{elli-con1}) and (\ref{elli-con2}) obviously. If we
can prove that $\Delta$ is weakly coercive, then all the
conclusions in section 2 hold for $\L=\Delta$ and its Green
function.\

Define the bilinear form $$
a_t(u,\varphi)=\int_\tm \langle \nabla
u,\nabla\varphi\rangle-\int t\langle u,\varphi\rangle
$$
from $H^{1,2}(\tm)\times H^{1,2}(\tm)$ to ${\mathbb C}$.
The form $a_t(u,\varphi)$ is bounded, since
$$
|a_t(u,\varphi)|\le ||u||_{H^{1,2}}\cdot ||\varphi||_{H^{1,2}},
$$
for $0\le t\le 1$. \

 The form $a_t(u,\varphi)$ is coercive, since
\begin{align*}
a_t(u,u)=&\int |\nabla u|^2-t\int |u|^2\ge \int |\nabla
u|^2-\frac{t}{\lambda_1-\delta}\int |\nabla u|^2\\
\ge &(1-\frac{t}{\lambda_1-\delta})\frac{1}{\lambda_1-\delta}\int
|u|^2,
\end{align*}
for $0<\delta<\lambda_1, 0\le t<\lambda_1-\delta$.

Hence
$$
a_t(u,u)\ge C_\delta ||u||^2_{H^{1,2}(\tm)}
$$
for $0\le t<\lambda_1-\delta$, where
$$
C_\delta=\frac{1}{2}(1-\frac{t}{\lambda_1-\delta})\min\{1,\frac{1}{\lambda_1-\delta}\}.
$$
If we take $\delta=\frac{\lambda_1}{2}$, then for any $0\le
t<\frac{\lambda_1}{3}$, there is
$$
a_t(u,u)\ge C_{\lambda_1}||u||^2_{H^{1,2}},
$$
where $C_{\lambda_1}=\frac{1}{6}\min\{1,\frac{2}{\lambda_1}\}$.
Now by Lax-Milgram theorem, for any $f\in H^{-1,2}(\tm)$, there
exists a unique $u\in H^{1,2}$ such that
$$
a_{\lambda_1/3}(u,v)=\langle f,v\rangle.
$$

Take $\varphi\ge 0, \varphi\in C^\infty_0(\tm)$, then the above
equality implies that
$$
(\Delta+\frac{\lambda_1}{3})u=-\varphi\le 0.
$$
On the other hand,
$$
a_{\lambda_1/3}(u^-,u^-)=a_{\lambda_1/3}(-u,u^-)=-\int \varphi
u^-\le 0.
$$
So by coercivity, there is $u^-=0$ and $u\ge 0$. Therefore if
$\varphi\neq 0$, we obtain a positive superharmonic function $u>0$
of the operator $\Delta+\lambda_1/3$.\

\begin{thm}\label{lapl-deca} There exist two positive numbers $C$ and $\alpha$
depending only on the geometry of $M$ such that $\forall (x,y)\in
\tm\times \tm$ and $d(x,y)\ge 1$, the following
holds:
\begin{equation} G(x,y)\le C e^{-\alpha d(x,y)}.
\end{equation}
\end{thm}
\begin{proof} This decay estimate was already proved in \cite{SY}. Here we give a different proof.
Firstly we prove that for $0<\epsilon<\lambda_1/3$,
and for any $x,y\in\tm$ satisfying $d(x,y)\ge 1$, we have
$G^\epsilon(x,y)\le C$, where $C$ only depends on $\lambda_1$.\

Let $f$ and $g$ be the characteristic function of the balls
$B(x,\rho)$ and $B(y,\rho)$ respectively, where
$\rho=\min\{r_0,1/3\}$. Then $G^\epsilon(fdv)$ is the solution of
the equation $\Delta u+\epsilon u=f$. By Schwarz inequality and
Lax-Milgram theorem, we have
$$
\int G^\epsilon(f)\cdot g\le (\int
|G^\epsilon(f)|^2)^{\frac{1}{2}}\cdot ||g||_{L^2}\le
C_{\lambda_1/3} ||f||_{L^2}||g||_{L^2}=C.
$$
Thus we have
$$
\int\int_{(\xi,\eta)\in B(x,\rho)\times
B(y,\rho)}G^\epsilon(\xi,\eta)d\xi d\eta\le C.
$$
Therefore there exists a point pair $(x_1,y_1)\in B(x,\rho)\times
B(y,\rho)$ such that
$$
G^\epsilon(x_1,y_1)\le C.
$$
Using Harnack inequality, we obtain
$$
G^\epsilon(x,y)\le C,
$$
for all $(x,y)$ such that $d(x,y)\ge 1$. Here $C$ only depends on
$M$. By Proposition \ref{decay-2}, we are done.
\end{proof}

\begin{crl} Given $\delta>0$, there exists $R=R(M,\delta)$ such
that $\forall x\in M$ and $r\ge R$, the $\Delta$-harmonic measure
$\mu_x$ of $x$ in $B(x,r)$ and the similar $\Delta+\epsilon$
harmonic measure $\mu^\epsilon_x$ satisfy
$$
\mu_x\le \delta\mu^\epsilon_x.
$$
\end{crl}
\begin{proof} It is a direct conclusion from corollary \ref{comp-meas} and
Theorem \ref{lapl-deca}.
\end{proof}

By Theorem \ref{lapl-deca}, the Green function $G$ of Laplace
operator vanishes at infinity. For the Green function of the
general elliptic operator $\L$, we need the following definition.
Let $\xi\in \tm(\infty)$, we say a function $u$ vanishes at $\xi$
in the $\L$-sense, if there exists a positive $\L$-superharmonic
function $w$ on $\tm$ such that $u=o(w)$ at $\xi$. If $\L(1)\le
0$, then the vanishing of $u$ at $\xi$ in the $\L$-sense is the
same as usual. It is shown in \cite{An}, page 509, that for any
$x\in\tm$, $G_x$ vanishes on $\tm$ in the $\L$-sense. Namely,
there exists a $\L$-superharmonic function $w$ such that
$G_x=o(w)$ at
 infinity.

\begin{prop}\label{poten-reduit} Let $0<\theta<\pi$, $\Gamma=\Gamma_{\sigma, t_0,
\theta},$ and $\Gamma_1=\Gamma_{\sigma, t_0+T_0, \theta}$ for some
$t_0$ and $T_0>0$. If $u(x)$ is a positive $\L$-harmonic function
in $\Gamma$ and vanishes in the $\L$-sense in $\tm(\infty)\cap
\Gamma$, then the reduit $u_1(x):=R^{\Gamma_1}_{u}(x)$ is a
$\L$-potential on $\Gamma$.
\end{prop}

\begin{proof} This is proved in \cite{An}, Theorem 2.
\end{proof}

\

\section{Hyperbolicity Estimates }
\setcounter{footnote}{0}

Let $\tm $ be a Cartan-Hadamard manifold with bounded geometry,
and recall the definition of $(h,T,\delta)$-non flat geodesics. We
have the following properties of the distance $\eta$:

\begin{prop}[Semi-additivity] \label {subad} For any $h>0$, we have
\begin{equation}
\eta_{\sigma([t_1,t_3])}(h)\ge
\eta_{\sigma([t_1,t_2])}(h)+\eta_{\sigma([t_2,t_3])}(h)
\end{equation}
for $t_1\le t_2\le t_3$.
\end{prop}

\begin{proof} Let $\phi:[t_1,t_3]\mapsto
\mathcal{U}_h$ be a path from $S^\perp_h(\sigma(t_1))$ to
$S^\perp_h(\sigma(t_3))$. For clear  topological reasons, the path
$\phi$ must intersect $S^\perp_h(\sigma(t_2))$ at $\phi(t^*)$.
Let $L(\phi|_{[s,s+\delta]})$ be the length of
$\phi|_{[s,s+\delta]}$. We have
$$
L(\phi|_{[t_1,t_3]})=L(\phi|_{[t_1,t^*]})+L(\phi|_{[t^*, t_3]})\ge
\eta_{\sigma([t_1,t_2])}(h)+\eta_{\sigma([t_2,t_3])}(h).
$$
\end{proof}

\begin{prop}[Continuity]\label {nonflat-cont} For fixed $t_1,t_2$ and $h$, the function $\eta_{\sigma([t_1,t_2])}(h)$ depends continuously on $\sigma'(0)$. Namely, for fixed $\delta_0$, there exists $\varepsilon  = \varepsilon (\tm, t_1,t_2,h, \delta_0)$ such that if $d_{S\tm}(\sigma'(0),\tau'(0)) < \varepsilon $, then  $$|\eta_{\sigma([t_1,t_2])}(h)-
\eta_{\tau([t_1,t_2])}(h)| \; < \; \delta_0.$$
\end{prop}

\begin{proof} Indeed if $\sigma '(0)$ and $\tau'(0) $ are close enough, then the closed sets
$$\mathcal{U}_h(\sigma([t_1,t_2])), \; S^\perp_h(\sigma(t_1))\; \text {and }\; S^\perp_h(\sigma(t_2)) $$ are sufficiently close to respectively the closed sets  $$\mathcal{U}_h(\tau([t_1,t_2])), \; S^\perp_h(\tau(t_1)) \; \text {and} \; S^\perp_h(\tau(t_2))$$  that the respective distances $$ d_{\mathcal{U}_h}(S^\perp_h(\sigma(t_1)),
S^\perp_h(\sigma(t_2)))\; \text { and} \; \; d_{\mathcal{U}_h}(S^\perp_h(\tau(t_1)),
S^\perp_h(\tau(t_2)))$$ are close. Moreover, by bounded geometry, if $t_1,t_2$ and $h$ are bounded, the explicit $\varepsilon $ of the above argument can be uniformly chosen, depending only on $\delta _0$.
\end{proof}

The other properties of ${\eta}$ we use need some explicitation:
 Let $F=\exp: \mathcal{N}(\sigma(\mathbb R))\mapsto \tm$ be
the exponential map (Fermi-map) along $\sigma$, where
$\mathcal{N}(\sigma(\mathbb R))$ is the normal bundle along $\sigma$.
If $\overrightarrow{Y}: {\mathbb R}\mapsto T_\sigma\tm$ is a $C^2$-smooth
vector field  along $\sigma$ with $\overrightarrow{Y}\perp\sigma'$
and $|\overrightarrow{Y}|\equiv1$, we consider the map
\begin{align*}
F=F_{\overrightarrow{Y}}: {\mathbb R}_+\times [t_1,t_2]&\mapsto
\widetilde{M}\\
(s,t)&\mapsto \exp_{\sigma(t)}(s\overrightarrow{Y}(t)).
\end{align*}
For fixed $h$, the map
$
F(s,t)=\exp_{\sigma(t)}[s\overrightarrow{Y}(t)], \forall (s,t)\in
[0,h]\times [t_1,t_2],
$
 gives a two-dimensional embedding surface with image
$\square_{t_1,t_2,h}$.
  Proposition 4.3 below implies that
 $F:{\mathbb R}^2\mapsto
\tm$ is a distance-increasing map, so the intrinsic curvature
$K_{\square_{t_1,t_2,h}}$ is well-defined. There is an intrinsic
curvature function
\begin{equation}\label{curv}
K_{\square_{t_1,t_2,h}}(s,t)=K(\frac{\partial F}{\partial
t},\frac{\partial F}{\partial s})-\frac{|\nabla_{\frac{\partial
F}{\partial s}}\frac{\partial F}{\partial t}|^2}{|\frac{\partial
F}{\partial s}\wedge \frac{\partial F}{\partial t}|^2}.
\end{equation}

This curvature function is related to the following length function:
$$
l(h)=L(F(h,\cdot)|_{[t_1,t_2]}).
$$
by the following proposition:
\begin{prop}\label{prop3.23} Let $\sigma, \overrightarrow{Y} $ and
$F=F_{\overrightarrow{Y}}$ be as above. Then
\begin{itemize}
\item[(1)] $l(h)$ is a convex function of $h$;

\item[(2)] If $r(x)=d(x,\sigma({\mathbb R}))$, then $\text{Hess}(r)(X, X) =
 \langle \nabla_X \nabla r, X \rangle \ge 0$ and
$$
\frac{dl}{dh}=\int^{t_2}_{t_1}\langle \nabla_{\frac{\partial
F}{\partial t}}\nabla r, \frac{\partial F}{\partial t}\rangle
\frac{1}{|\frac{\partial F}{\partial
t}|}dt=\int_{F(h,\cdot)}k_g(\cdot,h)dl \ge 0,
$$
where $k_g$ is the geodesic curvature of the curve $t\mapsto
F(h,t)$ with respect to $\nabla r$,
$k_g = - \langle \nabla_{\frac{\partial F}{\partial t}}(\frac{\partial F}{\partial t} ), \nabla r \rangle
= \text{Hess}(r)( \frac{\partial F}{\partial t}, \frac{\partial F}{\partial t}  )$.

\item[(3)]
$$
-\int_{\square_{t_1,t_2,h}}K_{\square_{t_1,t_2,h}}dA=\frac{\partial l}{\partial
h}.
$$

\item[(4)] $$ \frac{\partial l}{\partial h}(h)\ge
\frac{l(h)-l(0)}{h}$$

\item[(5)] $$ l(h) \; \ge \; t_2 - t_1. $$
\end{itemize}
\end{prop}

\begin{proof} Recall that $h\mapsto F(h,t) $ is a geodesic.
Therefore
$$
{J}^t(h)=\frac{\partial F}{\partial t}(h,t)
$$
is a Jacobi field along the geodesic ray $\Psi_t:h\mapsto
\Psi_t(h)=F(h,t)$.

It is easy to see that if $K_{\widetilde{M}}\le 0$, then the
function $h\mapsto ||J^t(h)||$ is a convex function in $h$, i.e.,
$$
\frac{\partial^2||J^t(h)||}{\partial h^2}\ge 0.
$$
Therefore
$$
\frac{\partial^2 l}{\partial
h^2}=\int^{t_2}_{t_1}\frac{\partial^2 ||\frac{\partial F}{\partial
t}||}{\partial h^2}dt\ge 0.
$$
For (2), it is a direct consequence of the first variational
formula, where $\nabla r|_{F(h,t)}=\frac{\partial F}{\partial
h}(h,t)$. In addition, it is proved in [BGS] that if $\sigma(\mathbb R)$ is a convex subset, then
$r(x)$ is a convex function in $x \in \tm$.

The assertion (3) follows from the Gauss-Bonnet formula on
$\square_{t_1,t_2,h}$. To see this we observe that
$||\overrightarrow{Y}(t)||=1$. It is clear that
$r(y)  \equiv h $ for all $y \in  {\mathcal{U}}_h(\sigma({\mathbb R}))$. It follows
that $r^{-1}(h) = {\mathcal{U}}_h(\sigma({\mathbb R}))$ and  $(\nabla
r|_{F(h,t)})\perp {\mathcal{U}}_h(\sigma({\mathbb R}))$.
Hence, we have a rectangle of curved top.

The discussion above implies that
\begin{equation}
\frac{\partial F}{\partial h}=\nabla r\perp \frac{\partial
F}{\partial t},
\end{equation}
because $\frac{\partial F}{\partial t}\in T_{F(t,h)}[\mathcal{U}_h(\sigma({\mathbb R}))]$.

Therefore,we apply the Gauss-Bonnet formula to get
$$
2\pi=\frac{\pi}{2}+\frac{\pi}{2}+\frac{\pi}{2}+\frac{\pi}{2}+\int_{F(h,\cdot)}k_g
dl+\int_{\square_{t_1,t_2,h}}K_{\square_{t_1,t_2,h}}dA.
$$
Thus,
$$
-\int_{\square_{t_1,t_2,h}}K_{\square_{t_1,t_2,h}}dA=\int_{F(h,\cdot)}k_g
dl=\frac{\partial l}{\partial h}.
$$
This proves (3).

For (4), we already proved that $l(h)$ is a convex function.
Thus, we have
$$
\frac{\partial l }{\partial h}\ge
\frac{l(h)-l(0)}{h}.
$$
Since, by (\ref {curv}), the Left Hand Side of (3) is nonnegative,
$\frac {\partial l}{\partial h} \ge 0$ and so $l(h) \ge l(0) = t_2 - t_1$. This proves (5).
\end{proof}

By definition we have: $$\eta_{\sigma ([t,t+T])} (h) :=\inf_{\stackrel{|\overrightarrow{Y}|=1}{\overrightarrow{Y}\perp\sigma'}}\{ L(F_{\overrightarrow{Y}}(h,\cdot)|_{[t,t+T]})\}.$$
Therefore,  by Proposition \ref {prop3.23}, if a geodesic $\sigma $ is $(h,T,\delta)$-non flat at $t$, then it satisfies:
\begin{equation}\label{non-flat}
\hat{K}_{\sigma,h}(t,
t+T):=\inf_{\stackrel{|\overrightarrow{Y}|=1}{\overrightarrow{Y}\perp\sigma'}}\{ \int\int_{\square_{t,t+T,h}}-K_{\square_{t,t+T,h}}(s,t)\left|\frac{\partial
F}{\partial s}\wedge \frac{\partial F}{\partial t}\right|ds\;dt\}\ge
\delta.
\end{equation}
We also have:
\begin{crl}\label{flat} If for some $t,T $ and $h$ a geodesic $\sigma $
satisfies $\eta_{\sigma ([t,t+T])} (h) = T$, then there is a field
$\overrightarrow {Y}$ along $\sigma $ such that the rectangle $\square_{t,t+T,h}$
 is totally geodesic and flat.
\end{crl}
\begin{proof}
This assertion was indeed implicitly stated in [BGS]. For the convenience of readers, we
present a short proof here. Let $\mathcal P_\sigma: \tm \mapsto \sigma(\mathbb R)$ be the
nearest point
projection. Since $\tm $ is a Cartan-Hadamard manifold and $\sigma(\mathbb R)$ is a closed
 convex subset,
it was proved in [BGS] that $\mathcal P_\sigma$ is a distance non-increasing map. Thus, we
have
$$
d_{\tm}(x, y) \ge d(\mathcal{P}_\sigma (x), \mathcal{P}_\sigma(y)    ).
$$
Equality holds in above inequality if and only if the four points
$\{ x, y, \mathcal{P}_\sigma (x), \mathcal{P}_\sigma(y)\}$
are vertices of a totally geodesic  flat rectangle $\Box$, see [BGS].

Suppose that $\eta_{\sigma ([t,t+T])} (h) = T$. By compactness, there is a point $ x \in S^\perp_h(\sigma(t))$,
a point $y \in   S^\perp_h(\sigma(t+T))$ and a shortest curve
on $ \mathcal{U}_h(\sigma([t,t+T]))$ realizing $d_{\mathcal{U}_h}(x,y)
 = \eta_{\sigma ([t,t+T])} (h) = T.$ Therefore, we have the following equalities and inequalities:
 $$
T = d_{\mathcal{U}_h}(x, y) \ge  d_{\tm}(x, y) \ge d(\mathcal{P}_\sigma (x), \mathcal{P}_\sigma(y)    ) = T.
 $$
Hence,  all inequalities above become equalities. In particular, we have  $
d_{\tm}(x, y) = d(\mathcal{P}_\sigma (x), \mathcal{P}_\sigma(y)    ),
$
which
implies that the four points
$\{ x, y, \mathcal{P}_\sigma (x), \mathcal{P}_\sigma(y)\}$
are vertices of a totally geodesic  flat rectangle $\Box_{t,t+T,h}$.
\end{proof}

We can describe  the geometric consequences of non-flatness we shall use.
Let $\widetilde{M}$ be a Cartan-Hadamard manifold and $\sigma$ be a
geodesic line of unit speed. Recall the family of cones
$
\Gamma_{\sigma,t,\theta}=\{x\in\widetilde{M}|\angle_{\sigma(t)}(\sigma'(t),x)<\theta\}.
$

\begin{prop}\label{cone1} Suppose that the geodesic $\sigma $ is $(h,T,  {\pi}/{4})$-non flat at 0.
Then:
\begin{equation}\label{cone-cont1}
\Gamma_{\sigma, T+h ,  {3\pi}/{4}}\subset
\Gamma_{\sigma, 0, {\pi}/{2}} \; \text {and}\;  \Gamma _{\sigma, T, {\pi}/{2}} \subset \Gamma_{\sigma, -h,  {\pi }/{4}}.
\end{equation}
\end{prop}

\begin{proof} Let us show the first inclusion, the proof of the other one is  similar. It suffices
to show that there is no geodesic triangle with  one side $\sigma
([0, T+h])$, another side $\tau $  in $\partial \Gamma_{\sigma,
T+h , {3\pi}/{4}}$ and the third side in $\partial \Gamma_{\sigma,
0, {\pi}/{2}}$. Suppose there is  such a rectangle
geodesic triangle $\triangle_{\sigma(0), \sigma(T+h), \tau(b)}$
with given three vertices $\{ \sigma(T+h),\sigma(0), \tau(b)  \}$,
where $\tau(b) \in [ (\partial \Gamma_{\sigma, T+h , {3\pi}/{4}})
\cap (   \partial \Gamma_{\sigma, 0, {\pi}/{2}}) ]$. We derive a contradiction as follows. We   choose  the vector field $\overrightarrow{Y}: {\mathbb
R}\mapsto  T_\sigma\tm$  along $\sigma$ with
$\overrightarrow{Y}\perp\sigma'$ and $|\overrightarrow{Y}|\equiv1$
in such a way that $\exp_{\sigma (t)} [S(t)\overrightarrow{Y}] $
lies in $\tau $ for some $S(t)$. As before, we let $\square
_{0,T,h} = \{ \exp_{\sigma (t)} [s \overrightarrow{Y}]| 0\le t\le
T, 0\le s \le h  \}$. By comparison with the Euclidean plane, we
have $S(t) \ge h$ for $0 \le t \le T$. Therefore the triangle
$\triangle_{\sigma(0), \sigma(T+h), \tau(b)} = \{ \exp_{\sigma
(t)} [s \overrightarrow{Y}]| 0\le t\le T+h, 0\le s \le S(t) \}$
contains the subset $\square _{0,T,h}$. By Proposition 4.3
(3)-(4), we  have
$$
-\int_{\triangle_{\sigma(0), \sigma(T+h), \tau(b)}
}K_{\triangle}dA \ge
-\int_{\square_{t_1,t_2,h}}K_{\square_{t_1,t_2,h}}dA=\int_{F(h,\cdot)}k_g
dl=\frac{\partial l}{\partial h} > \frac{\pi}{4}.
$$
This together with the Gauss-Bonnet formula implies that the
sum of inner angles  of $\triangle_{\sigma(0), \sigma(T+h),
\tau(b)}$
 is smaller than $ (\pi -  \frac{\pi}{4}) = {3\pi }/{4}$, which is impossible.
\end{proof}

The same proof also yields:
\begin{prop}\label{cone2} Let $\varepsilon  >0$, and suppose that there is $t_+, t_- >
\frac{h}{\tan \varepsilon}  $  such that the  geodesic $\sigma $ is $(h,T,  {\pi}/{2})$-non flat at $t_+$ and $-t_-$. Then:
\begin{equation}\label{cone-cont2}
\Gamma_{\sigma, T+t_+ ,  {\pi}/{2}}\subset
\Gamma_{\sigma, 0,\varepsilon } \; \text {and } \; \Gamma_{-\sigma, t_- ,  {\pi}/{2}}\subset
\Gamma_{-\sigma, -T ,\varepsilon }.
\end{equation}
\end{prop}

\

The main geometric estimate related to the Martin boundary is Ancona's $\Phi$-chain condition. For a cone $\Gamma _{\sigma, 0, \theta}$, it says that one can find a time $T_0$  such that,  for $x\in \partial \Gamma_{\sigma, 0,\theta}$,
\begin{equation}\label{prob-ineq1}
d(x, \Gamma_{\sigma ,T_0,\theta}) \to \infty \; {\text as} \; d(x,\sigma (0)) \to \infty.
\end{equation}

When $\widetilde{M}={\mathbb R}^n$ is the Euclidean space, then for $x \in \partial \Gamma_{\sigma,0,\theta}$, $
d(x,
\Gamma_{\sigma,T_0,\theta})\le T_0.
$ and
 can NOT be unbounded. For the same reason, if $\sigma({\mathbb R})$ is a boundary
of a totally geodesic flat half plane ${\mathbb R}^2_+$, then
(\ref{prob-ineq1}) fails on $ {\mathbb R}^2_+\cap \partial
\Gamma_{\sigma,0,\theta}$. However, the cone property
(\ref{cone-cont1}) implies  a stronger form of (\ref{prob-ineq1}).

\begin{prop}\label{fichain} Let $0<\theta\le \frac{\pi}{2}$. If $\Gamma_{\sigma, T_0, \theta+\varepsilon'_0}\subset
\Gamma_{\sigma, 0 ,\theta}$ for some $T_0>0$ and $\varepsilon'_0>0$,
then
\begin{equation}
d(x,\Gamma_{\sigma, T_0,\theta})\ge
\varepsilon_0[d(x,\sigma(0))-\frac{1}{\varepsilon_0}]
\end{equation}
for $x\in \partial \Gamma_{\sigma,0,\theta}$ and  some
$\varepsilon_0>0$ which depends only on $\varepsilon'_0$ and $T_0$.

\end{prop}

\begin{proof} By our assumption, if $x\in \partial
\Gamma_{\sigma,0,\theta}$ then
\begin{equation}\label{angl-larg1}
\angle_{\sigma (T_0)}(x,\Gamma_{\sigma,T_0,\theta})\ge\varepsilon'_0.
\end{equation}
Recall that $\exp_{\sigma(T_0)}: {\mathbb R}^n\mapsto \widetilde{M}$
is a distance increasing map. If
$$
\stackrel{\circ}{\Gamma}\!_{\sigma,T_0,\theta}=\{\overrightarrow{\omega}\in
T_{\sigma(T_0)}\widetilde{M}|\angle(\overrightarrow{\omega},
\sigma'(T_0))\le \theta\}
$$
and $\overrightarrow{u}_x=\exp^{-1}_{\sigma(T_0)}x$, then by
(\ref{angl-larg1}) we have
\begin{equation}
d_{{\mathbb R}^n}(\overrightarrow{u}_x,
\stackrel{\circ}{\Gamma}\!_{\sigma,T_0,\theta})\ge
|\overrightarrow{u}_x|\sin \varepsilon'_0.
\end{equation}
Since $\exp_{\sigma(T_0)}:{\mathbb R}^n\mapsto \widetilde{M}$ is
distance increasing, we conclude that
$$
d_{\widetilde{M}}(x, {\Gamma}\!_{\sigma, T_0,\theta})\ge
d(x,\sigma(T_0))\sin \varepsilon'_0\ge [d(x,\sigma(0))-T_0]\sin
\varepsilon'_0.
$$
Then we choose
$\varepsilon_0=\min\{\sin\varepsilon'_0,\frac{1}{T_0\sin\varepsilon'_0}\}$ and we obtain (4.8).
\end{proof}

\

\section{Boundary Harnack Inequality and Martin boundary}\
\setcounter {footnote}{0}

\subsection{Boundary Harnack Inequality, proof of Theorem \ref{barrier}}\

We assume in this section that the geodesic $\tau $ admits a $(h,T,R) $ barrier, and we are going to prove (\ref {unif-G}). The proof of (\ref {unif-cons}) is the same.

 \begin{prop}\label{phi-cond}
Assume the geodesic $\tau : {\mathbb R} \mapsto \tm$ is $(h,T,
\pi/4)$-non flat at 0 and set $T_0 = T+h$.  Denote $ x_p =\tau(p
T_0),p\ge 0$. Then there exists a constant $C=C(\tm, h, T )$ such
that
\begin{equation}\label{proof-113}
G(y,x_p)\le C G(x_0,x_p)G^\epsilon(y,x_1), $$
$$
\forall y\in \tm\setminus \Gamma_{\tau, 0 ,\pi/2},\;
\forall p\ge 1.
\end{equation}
Furthermore, for any $x\in \overline{x_px_{p+1}}$, the line
segment between $x_p$ and $x_{p+1}$($p\ge 1$), one has
\begin{equation}
G(y,x)\le C G(x_0,x)G^\epsilon(y,x_1), \forall y\in \tm\setminus
\Gamma_{ \tau, 0 ,\pi/2}.
\end{equation}
\end{prop}

\begin{proof}
We denote  $\Gamma=\Gamma_{\tau,0,\pi/2},
\Gamma_1=\Gamma_{\tau,T_0, \pi/2}$. \

To prove (\ref{proof-113}), we firstly prove the following
inequality: there exists $C_p=C(p, \tm, h,T)$ such that
\begin{equation}\label{proof-114}
G(y,x_p)\le C_p G(x_0,x_p) G^\epsilon(y,x_1),\forall y\in
\tm\setminus \Gamma.
\end{equation}
By Remark \ref{deca-rema2}, we have, with $C''_p = ce^{-\beta p
T_0}$ for some $c,\beta$ depending only on $(\tm,\L)$ and
\begin{equation}\label{proof-115}
G(x_0,x_p)\ge C''_p.
\end{equation}
By construction,
$B(x_1,h)\subset \Gamma$.

Take $y_0\in \partial B(x_1,h)$, by Harnack
inequality, then we obtain
$$
G(y_0,x_p)\le C'_p G(y_0,x_1).
$$
Applying the Harnack inequality to the variable $y$, we have
\begin{equation}\label{proof-116}
G(y,x_p)\le C'_{p,1}G(y,x_1), \forall y\in \partial
B(x_1,h),
\end{equation}
with $C'_{p,1} = C'_{p,1} (\tm, d(x_1,x_p))$.
Similarly we can prove that
\begin{equation}\label{proof-117}
G(y,x_p)\le C'_{p,2}G(y,x_1), \forall y\in \partial
B(x_p,h),
\end{equation}
with $C'_{p,2} = C'_{p,2} (\tm, d(x_1,x_p))$.
Let $C_p=\max\{C'_{p,1},C'_{p,2}\}$. Combining (\ref{proof-116})
and (\ref{proof-117}), there is
$$
G(y,x_p)\le C_p G(y,x_1), \forall y\in \partial
B(x_1,h)\cup \partial
B(x_p,h).
$$
Now using maximum principle, we have
$$
G(y,x_p)\le C_p G(y,x_1), \forall y\in \tm\setminus
B(x_1,h)\cup B(x_p,h).
$$
In particular, we have
\begin{equation}\label{proof-118}
G(y,x_p)\le C_p G(y,x_1)\le C_p G^\epsilon(y,x_1), \forall y\in
\tm \setminus \Gamma.
\end{equation}
By (\ref{proof-115}) and (\ref{proof-118}), we obtain
(\ref{proof-114}). \

The  proof of Proposition \ref {phi-cond} will  consist in showing
that one can take the constant in (\ref{proof-114}) independent of
$p$.  Observe indeed  that to obtain (\ref{proof-114}), we only
used the relative distances of $x_0, x_1 $ and $ x_p$ and that
$B(x_1,h)\subset \Gamma$.  Therefore (\ref{proof-114}) can be
applied to the cone $\Gamma_1$ to get
\begin{equation}\label{proof-119}
G(y,x_{p+1})\le C_p G(x_1,x_{p+1}) G^\epsilon(y,x_2), \forall y\in
\tm \setminus \Gamma_1.
\end{equation}

Applying the Harnack inequality of $\L^*+\epsilon I$ to its Green
function $G^\epsilon(y,x)$, one has
$$
G^\epsilon(y,x_2)\le C' G^\epsilon(y,x_1),\forall y\in \partial
B(x_2, h).
$$
Then maximum principle and $B(x_2,h)\subset \Gamma _1$   implies that
\begin{equation}\label{proof-120}
G^\epsilon(y,x_2)\le C' G^\epsilon(y,x_1), \forall y\in
\tm\setminus \Gamma_1,
\end{equation}
where $C'=C'(\tm,T+h)$. By Harnack inequality, we have
\begin{equation}\label{proof-121}
G(x_1,x_{p+1})\le C'' G(x_0,x_{p+1}),
\end{equation}
where $C''=C''(\tm,T+h)$ is independent of $p$ for $p\ge 1$.

Combining (\ref{proof-119}), (\ref{proof-120}) and
(\ref{proof-121}), one has
\begin{align}\label{proof-122}
G(y,x_{p+1})\le& C_p C'' G(x_0,x_{p+1}) C' G^\epsilon(y,x_1)\nonumber\\
=&C_p C G(x_0,x_{p+1}) G^\epsilon(y,x_1),\forall y\in \tm\setminus
\Gamma_1,
\end{align}
where $C=C(\tm,h,T)$.

Now fixing $\displaystyle \delta=\frac{1}{C}$, by Corollary
\ref{comp-meas}, $\forall \epsilon>0$, there exists
$R_0=R_0(\tm,h,T,\epsilon)$ such that $\forall x\in \tm$ ,and
$r\ge R_0$, the harmonic measures $\mu_x$ and $\mu_x^\epsilon$ on
balls of radius $r$ about $x$ satisfy
\begin{equation}\label{proof-123}
\mu_x\le \delta \mu_x^\epsilon.
\end{equation}

Since the geodesic $\tau $ is $(h,T, \pi/4)$-non flat at 0, by
Propositions \ref {cone1} and \ref {fichain} there is a
$\varepsilon_0$ depending only on $T,h$ such that
\begin{equation}
d(y,\partial\Gamma_1)\ge
\varepsilon_0(d(y,x_0)-\frac{1}{\varepsilon_0}).
\end{equation}
Now we can take $\rho_1=\rho_1(\varepsilon_0,
R_0)=\rho_1(\tm,h,T,\epsilon)$ such that for any $y\in
\tm\setminus \Gamma$ and $d(y,x_0)\ge \rho_1$ the following holds:
\begin{equation}\label{proof-124}
d(y,\partial\Gamma_1)\ge R_0.
\end{equation}
For such $y$, the ball $B(y,R_0)\subset \tm\setminus \Gamma_1$. By
(\ref{proof-122}), for any $z\in \partial B(y,R_0)$, we have
$$
G(z,x_{p+1})\le C_p C G(x_0,x_{p+1}) G^\epsilon(z,x_1).
$$
So
\begin{align*}
&\int_{\partial B(y,R_0)} G(z, x_{p+1})d\mu_y(z)\\
&\le C_p C\int_{\partial B(y,R_0)} G^\epsilon(z,x_1) d\mu_y(z)
G(x_0,x_{p+1})\\
&\le C_p C \delta\int_{\partial B(y,R_0)} G^\epsilon(z,x_1)
d\mu_y^\epsilon(z) G(x_0,x_{p+1})\\
&\le C_p G(x_0,x_{p+1})G^\epsilon(y,x_1),
\end{align*}
i.e.,
\begin{equation}\label{proof-125}
G(y,x_{p+1})\le C_p G(x_0,x_{p+1}) G^\epsilon(y,x_1),
\end{equation}
for any $y\in \tm\setminus \Gamma$ and $d(y,x_0)\ge \rho_1$.

There exists $C=C(\tm, T_0)$ such that $G(x_0,x_1) C \ge 1$. Thus
\begin{equation}\label{proof-126}
G(y,x_{p+1})\le C G(x_0,x_{p+1}) G(y,x_1)
\end{equation}
at $y=x_0$.
Using Harnack inequality in the compact set $(\tm\setminus
\Gamma)\cap B_{\rho_1}(x_0)$, one has
\begin{equation}\label{proof-127}
G(y,x_{p+1})\le C' G(x_0,x_{p+1}) G(y,x_1),
\end{equation}
for any $y\in (\tm\setminus \Gamma)\cap B_{\rho_1}(x_0)$, where
$C'$ depends on $\rho_1$ and $C=C(\tm, h,T)$ in (\ref{proof-126}),
and hence depends only on $\tm, h, T$, but not on the
 the integer $p$.

Combining (\ref{proof-127}) and (\ref{proof-125}), one obtains
\begin{equation}\label{proof-128}
G(y,x_{p+1})\le \max\{C_p,C'\} G(x_0,x_{p+1}) G^\epsilon(y,x_1).
\end{equation}
So we can improve $C_{p+1}$ such that $C_{p+1}=\max\{C_p,C'\}$.

Hence we can take a uniform constant $C=\max\{C_1, C'\}$ such that
for any $p\ge 0$ the following inequality holds:
$$
G(y,x_{p+1})\le C G(x_0,x_{p+1}) G^\epsilon (y,x_1), \forall y\in
\tm\setminus \Gamma\;\text{and}\;\forall p\ge 0,
$$
where $C$ depends only on $\tm ,h, T$.

For  the general point $x\in \overline{x_px_{p+1}}$, one can use
Harnack inequality because of $ d(x_p,x_{p+1}) = T_0 $.
\end{proof}

Recall from Section 2 the definition of $\varepsilon^* = \varepsilon ^* (\tm, h, T)$.

\begin{crl}\label{phi-gree} Assume the geodesic $\sigma $ is $(h,T,\pi/2)$-non flat at 0. Set $\Gamma = \Gamma _{\sigma, 0, \pi/2 + \varepsilon ^*}$ and $\Gamma _1 =  \Gamma _{\sigma, 0, + \varepsilon ^*}\cap \Gamma _{\sigma , h+T, \pi/2}$. Then,  for any
$y\in \tm\setminus \Gamma$, for any $x\in \Gamma_1$, there is
$$
G(y,x)\le C G(x_0,x) G^{\epsilon}(y,x_1),
$$
where $C=C(\tm , h, T), x_0=\sigma(0)$, and
$x_1=\sigma(h+T)$.
\end{crl}

\begin{proof} By our choice of $\varepsilon ^*$ and Proposition \ref{phi-cond},  any geodesic  $\tau$ which satisfies $$ \tau (0) \;= \; \sigma (0) \; \text{and} \; \angle_{\sigma(0)}(\tau'(0),\sigma'(0)) \; <\; \varepsilon ^*,$$ is $(h,T,\pi/4)$-non flat at 0.
By  Proposition
\ref{phi-cond} for any
$x\in \tau ([h+T,+\infty))$ and any $y\in
\tm\setminus\Gamma_{\tau,0 ,\pi/2 }$, there is
\begin{equation}\label{proof-crl4.7-1}
G(y,x)\le C G(x_0,x) G^{\epsilon}(y,\tau (h+T)),
\end{equation}
where $C=C(\tm ,h,T)$. On the other hand, by comparison   $d(x_1,
\tau (h+T)) \leq \varepsilon ^* \sinh (K(h+T))$, where $-K$ is a
lower bound for the sectional curvature on $\tm$,  so that  by
Harnack inequality, there is a $C=C(\tm, h,T)$ such that for any
$y\in\tm\setminus\Gamma_{\tau,0, \pi/2}$ the following holds
\begin{equation}\label{proof-crl4.7-2}
G^\epsilon (y,\tau (h+T))\le C  G^{\epsilon}(y,x_1).
\end{equation}

Combining (\ref{proof-crl4.7-1})  and
(\ref{proof-crl4.7-2}), we have the conclusion for any $x\in \Gamma _{\sigma , 0, \varepsilon^*}$ at distance at least $h+T$ from $\sigma (0)$  and any $y\in
\tm\setminus \cup_\tau \Gamma _{\tau, 0, \pi/2} $,  in particular for points $x\in \Gamma _1$ and $y\in \tm \setminus \Gamma $.
\end{proof}

By Remark \ref{reverse}, we can apply  Corollary \ref {phi-gree}  to $-\sigma $ and  get:

\begin{crl}\label{gree-phi} Assume the geodesic $\sigma $ is $(h,T,\pi/2)$-non flat at 0. Set
 $\Gamma ' = \Gamma _{-\sigma, -T, \pi/2 + \varepsilon ^*}$ and $\Gamma '_1 =  \Gamma _{-\sigma, -T, + \varepsilon ^*}\cap \Gamma _{-\sigma , h, \pi/2}$. Then,  for any
$y\in \tm\setminus \Gamma '$, for any $x\in \Gamma'_1$, there is
$$
G(y,x)\le C G(x_0,x) G^{\epsilon}(y,x_1),
$$
where $C=C(\tm , h, T), x_0=\sigma(T)$, and
$x_1=\sigma(-h)$.
\end{crl}

We can now prove Theorem \ref {barrier}:

\begin{proof} Since the geodesic $\sigma $ admits a $(h,T,R)$ barrier, there is $t_4, t_5$, with
$h/\tan \varepsilon ^* \le t_5-t_4 \le T_1+R$ such that $\sigma $
is $(h,T, \pi/2)$-non flat at $t_4$. By Proposition \ref {cone2},
then
$$ \Gamma _{-\sigma , -t_4, \pi/2 } \; \subset \; \Gamma _{-\sigma , -t_5 -T, \varepsilon ^*}.$$
Moreover, there is $t_6$, with $ T+\frac{h}{\tan\epsilon^*}  \le
t_6 - t_5 \le R$ such that $\sigma $ is $(h,T, \pi/2)$-non flat at
$t_6$ and by Proposition \ref {cone2},
$$ \Gamma _{\sigma , t_6 +T, \pi/2 } \; \subset \; \Gamma _{\sigma , t_5 +T, \varepsilon ^*}.$$

Applying Corollary  \ref{gree-phi} we get for any  $y\in  \Gamma _{\sigma , t_6 +T, \pi/2 }$
and for any $x\in \tm \setminus  \Gamma _{\sigma , t_4, \pi/2 }$, there is
\begin{equation}
G(y,x)\le C_1 G(\sigma (t_5+T),x) G^{\epsilon}(y,\sigma (t_5-h)),
\end{equation}
where $C_1=C_1(\tm ,h,T)$. Using Harnack inequality, we have with a different  $C_1=C_1(\tm ,h,T)$:
\begin{equation}\label{ineq-otherside1}
G(y,x)\le C_1 G(\sigma (t_6+T),x) G^{\epsilon}(y,\sigma (0)),
\end{equation}

In the same way, using that the geodesic $\sigma $ is $(h,T,\pi/2)$-non flat at $t_1, t_2$ and $t_3$, and Corollary \ref {phi-gree},
we can obtain that for any $x\in \tm\setminus
\Gamma_{\sigma , t_1, \pi/2}$ and for any $y\in
\Gamma_{\sigma, t_3 +T, \pi/2}$, there is
\begin{equation}
G(x,y)\le C_2 G^\epsilon(x,\sigma (t_2 +T+h)) G(\sigma (t_2),y),
\end{equation}
and
\begin{equation}\label{ineq-oneside2}
G(x,y)\le C_2 G^\epsilon(x,\sigma (0)) G(\sigma (t_1),y),
\end{equation}where $C_2=C_2(\tm , h,T)$.

Set
$x_0 =\sigma (t_1),  x' = \sigma (0) $ and $ x_1 =  \sigma (t_6 +T)$,
and let $\Gamma=\Gamma_{\sigma ,t_1
{\pi}/{2}}$, $\Gamma' =  \Gamma_{\sigma ,0,{\pi}/{2}}$ and $\Gamma _1 =
\Gamma_{\sigma , t_6 +T, {\pi}/{2}}$. We claim that
the cone pair $\{\Gamma,\Gamma_1\}$ satisfies the conclusion of Theorem \ref {barrier}. Since $-T_2 \le  t_1$ and $t_6 +T \le  T_2$, Theorem \ref{barrier} will follow.

\

We follow \cite {An}. Let $y\in \tm\setminus\Gamma$ and $x\in \Gamma_1$. We have the
representation
$$
\Gstar_y(x)=G(y,x)=\int_{\partial \Gamma'}G(y,z)d\mu_x(z),
$$
where $\mu_x$ is a positive measure supported on $\partial\Gamma'$
and such that
$$
G(\mu_x)=R^{\partial\Gamma'}_{G_x}=R^\omega_{G_x},
\;\omega=\tm\setminus\bar{\Gamma}'.
$$
Now applying inequality (\ref{ineq-otherside1}) to $\Gstar$, we
have
\begin{align}\label{proof-129}
G(y,x)\le&
C\int_{\partial\Gamma'}G(y,x')G^\epsilon(x_0,z)d\mu_x(z)\nonumber\\
\le&
CG(y,x_1)\int_{\partial\Gamma'}\Gstar{}^\epsilon_{x_0}(z)d\mu_x(z).
\end{align}

Since $G^\epsilon(x_0,z)=\Gstar{}^\epsilon_{x_0}(z)$ is an
$\L^*$-potential, so $\Gstar{}^\epsilon_{x_0}(z)=\Gstar(\lambda)$
for some positive measure $\lambda$ on $\tm$. We have
\begin{equation}\label{proof-130}
\int_{\partial\Gamma'}\Gstar{}^\epsilon_{x_0}(z)d\mu_x(z)
=\int_\tm\int_{\partial\Gamma'}G(y,z)d\mu_x(z)d\lambda(y)=
\int_\tm G(\mu_x)d\lambda(y)
\end{equation}
By inequality (\ref{ineq-oneside2}), for any $z\in
\partial\Gamma'$, there is
\begin{equation}\label{proof-131}
G_x(z)\le C G_x(x')G^\epsilon_{x_1}(z)\le
CG_x(x_0)G^\epsilon_{x_1}(z).
\end{equation}
By the definition of reduit, we have by (\ref{proof-131}),
\begin{equation}\label{proof-132}
G(\mu_x)=R^{\partial\Gamma'}_{G_x}\le G_x\le
CG(x_0,x)G^\epsilon_{x_1}.
\end{equation}
Substitute (\ref{proof-132}) into (\ref{proof-130}) and then plug
in (\ref{proof-129}), then we can obtain
\begin{equation}\label{proof-133}
G(y,x)\le CG(y,x_1)G(x_0,x)\int_\tm G^\epsilon_{x_1}d\lambda.
\end{equation}
Now we only need to show that $\int_\tm G^\epsilon_{x_1}d\lambda$
is bounded. By Lemma \ref{anco-lm1}, there is
$$
\Gstar{}^\epsilon_{x_0}=\Gstar{}_{x_0}+\epsilon
\Gstar(\Gstar{}^\epsilon_{x_0}).
$$
Since $\Gstar{}^\epsilon_{x_0}=\Gstar(\lambda)$, so
$$
\lambda=-\L^*(\Gstar{}^\epsilon_{x_0})=\delta_{x_0}+\epsilon\Gstar{}^\epsilon_{x_0}dv_\tm.
$$
Thus
\begin{align}\label{proof-134}
\int\Gstar{}^\epsilon_{x_1}d\lambda(z)=&G^\epsilon(x_1,x_0)+\epsilon\int
G^\epsilon(z,x_1)G^\epsilon(x_0,z)dv_\tm(z)\nonumber\\
=&G^\epsilon(x_1,x_0)+\epsilon
\Gstar{}^\epsilon(\Gstar{}^\epsilon_{x_1})(x_0)\nonumber\\
=&\Gstar{}^\epsilon_{x_1}(x_0)+\epsilon
\Gstar{}^\epsilon(\Gstar{}^\epsilon_{x_1})(x_0)
\end{align}
Using Lemma \ref{anco-lm1} again to the operator $\L^*+\epsilon I$
and $\frac{\epsilon}{2}$ instead of $\L$ and $\epsilon$, we have
\begin{equation}\label{proof-135}
\Gstar{}^\epsilon_{x_1}(x_0)+\epsilon\Gstar{}^\epsilon
(\Gstar{}^\epsilon_{x_1})(x_0)\le \widehat{G}(x_0,x_1)\le C,
\end{equation}
where $\widehat{G}(x,y)$ is the Green function of
$\L^*+\frac{\epsilon}{2}$. By Proposition \ref{anco-prop-7}, the
constant $C$ here only depends on $M,\L$ and the distance between
$x_0$ and $x_1$ and hence depend only on $M, \L,\theta, \tau$.\

Combining (\ref{proof-133}), (\ref{proof-134}) and
(\ref{proof-135}), the whole proof is finished.
\end{proof}

\subsection{Proof of Proposition \ref{unicity}}\

\begin{proof} Let $u(x),v(x)\in C_\xi$. Since $v(x)=O(G_{x_0})=o(w),$ for some $\L$-superharmonic function,
 for $x\in \bar{\Gamma}_{-
{\tau},-t_k-2T_2,\pi/2}\cap \tm (\infty)$, the reduit of $v_1(x)$
of $v(x)$ on $\Gamma_{- {\tau},-t_k-T_2,\pi/2}$ with respect to
$\Gamma_{- {\tau},-t_k-2T_2,\pi/2}$ is a potential by Proposition
\ref{poten-reduit}. So it has the representation formula:
$$
v_1(x)=\int_{\partial\Gamma_{- {\tau}, -t_k-T_2,\pi/2}}
g^k(x,y)d\nu^k(y), \forall x\in \Gamma_{- {\tau},-t_k +T_2,\pi/2},
$$
where $\nu_k(y)$ is positive measure on $\partial\Gamma_{-
{\tau},-t_k-T_2 ,\pi/2}$. According to the definition of reduit,
we have
\begin{equation}\label{redu-repr}
v(x)=v_1(x)\le C\int_{\partial\Gamma_{- {\tau},t_k, \pi/2}}
g^k(x,\tau(t_k))g^k(\tau(t_k-T_2),y)d\nu^k(y), \forall x\in
\Gamma_{- {\tau},-t_k+T_2,\pi/2},
\end{equation}
where we used (\ref{unif-cons}) and the constant $C$ here is
independent of $k$.

By (\ref{redu-repr}), there is
\begin{equation}\label{harn-v}
v(x)\le C g^k(x,\tau(t_k+T_2)) v(\tau(t_k-T_2)).
\end{equation}
On the other hand, by Harnack inequality and the  maximum principle, we have
$$
u(x)\ge C' u(\tau(t_k-T_2)) g^k(x,\tau(t_k+T_2)), \forall x\in
\Gamma_{- {\tau},-t_k-2T_2,\pi/2} \setminus B(\tau(t_k+T_2),\frac{T_2}{2}).
$$
In particular, one has
\begin{equation}\label{harn-u}
u(x)\ge C' u(\tau(t_k -T_2))g^k(x,\tau(t_k+T_2)), \forall x\in
\Gamma_{- {\tau},-t_k+T_2,\pi/2}.
\end{equation}
Combining (\ref{harn-v}) and (\ref{harn-u}), there is
\begin{equation}\label{rate-uv1}
\frac{v(x)}{u(x)}\le C\frac{v(\tau(t_k -T_2))}{u(\tau(t_k-T_2))},
\end{equation}
where $C$ is independent of $k$. Similarly, one has
\begin{equation}\label{rate-uv2}
C^{-1}\frac{v(\tau(t_k-T_2))}{u(\tau(t_k-T_2))}\le \frac{v(x)}{u(x)}\le
C\frac{v(\tau(t_k-T_2))}{u(\tau(t_k-T_2))}, \forall x\in \Gamma_{-
{\tau}, -t_k+T_2,\pi/2}.
\end{equation}
Let $x=\tau(t_0)$, then
$$
\frac{v(x)}{u(x)}\le
C^2\frac{v(\tau(t_0))}{u((\tau(t_0)))}:=\lambda,\forall x\in
\Gamma_{- {\tau},-t_k +T_2,\pi/2}.
$$
Since $\lambda$ is independent of $k$, we can let $k\to
\infty$ and obtain
$$
\frac{v(x)}{u(x)}\le \lambda,\forall x\in\tm,
$$
which implies $\dim C_\xi\le 1$.
\end{proof}

\subsection{Proof of Proposition \ref{existence}}\

\begin{proof}
 Let
$x_k=\tau(t_k)$, where $t_k \to \infty$ are the  barrier times. For $z_i \to \xi $ as $i\to \infty $,  define
$k^i_{x_q}:=\frac{G(x,z_{i})}{G(x_q,z_{i})}$. Then
$k^i_{x_q}(x)$ is a sequence of positive harmonic functions on
$\tm$ and that satisfy  the normalization condition at $x_q:\;
k^i_{x_q}(x_q)=1$. Therefore there exists a subsequence such that
$\lim_{i\to\infty} k^i_{x_q}(x)=k_{x_q}(x)$ exists (where we still
use the same index $i$).

 For $i$ large enough, $z_i \in \Gamma _{\tau, t_q +T_2, \pi/2}$ and, setting $\Gamma _q = \Gamma _{\tau, t_q -T_2, \pi/2}$,  Theorem \ref{barrier} gives that
\begin{equation}\label{exis-proo-2}
k_{x_q}(x)\le CG(x,x_{q}), \forall x\in \tm\setminus \Gamma_q.
\end{equation}

Equation (\ref{exis-proo-2}) means that
\begin{equation}\label{exis-proo-3}
k_{x_q}(\xi')=0 \;\text{ in the $\L$-sense},\;\forall \xi'\in
\tm(\infty)\setminus \overline{\Gamma_q}.
\end{equation}
Similarly we can consider the positive harmonic function
normalized at $x_0$, i.e., $k_{x_0}(x)$.

By Harnack inequality, there exists a constant $C_q$ such that for
sufficiently large $i$ the following holds:
$$
C_q^{-1}G(x_q,z_{i})\le G(x_0,z_{i})\le C_q G(x_q, z_{i}).
$$
Therefore we have
$$
\frac{G(x,z_{i})}{G(x_0,z_{i})}\le C_q
\frac{G(x,z_{i})}{G(x_q,z_{i})},
$$
for large $i$.

So we have
\begin{equation}\label{exis-proo-5}
k_{x_0}(x)\le C_q k_{x_q}(x)\le C_q G(x,x_{q}),\forall x\in
\tm\setminus\Gamma_q.
\end{equation}
On the other hand, by Proposition \ref{cone2}, for any $\varepsilon>0$, there is
a $q\in {\mathbb N}$ such that
$\Gamma_q\subset\Gamma_{\tau, 0,\varepsilon}$. Combining this
fact and (\ref{exis-proo-5}), we obtain
$$
k_{x_0}(\xi')=0  \;\text{ in the $\L$-sense},\;\;\forall \xi'\in \tm (\infty)\setminus \{\xi\}.
$$
Since $k_{x_0}(x_0)=1$, $k_{x_0}$ is a nontrivial positive
harmonic function on $\tm$, and hence is a Poisson kernel
function at $\xi$.
\end{proof}

\

\section{Abundance of Martin points}\

We assume in this section that $\tm $ is the universal cover of a compact Riemannian manifold
$M$ of class $C^3$, nonpositive  curvature and geodesic rank one.
For each $v \in T\tm$, we can write $v = (x, \vec \theta)$ with
$\vec \theta \in T_x \tm$. If $Exp: T\tm \mapsto \tm$ is the
exponential map given by $Exp(v) = Exp_x(\vec \theta)$. We always
write
$$
\sigma _v (t) = Exp_x ( t \vec \theta) $$
and
$$ \varphi_s (v) = (Exp_x ( s \vec \theta),
\frac{d[Exp_x ( s \vec \theta) ]}{ds}) = (\sigma _v (s), \sigma_v'(s)).
$$
The 1-parameter family of diffeomorphisms $\{ \varphi_s \}$ is called the geodesic flow on $S\tm$.

Let $\pi: S\tm \mapsto \tm$ be the foot-point projection. For any given $v
\in S\tm$, we consider the Busemann function
$$
   \hat b_v(x) = \lim_{t \to + \infty} [ d(\sigma _v (0), \sigma _v (t)) - d(x, \sigma _v (t))]
       = \lim_{t \to + \infty} [ t  - d(x, \sigma _v (t))],
$$
for $x \in \tm$.

The the level set $\Sigma_v = \hat b_v^{-1}(0)$ is called a horosphere with the inner normal vector
$v$. We also let
$$
\mathcal{H}_v = \{(y, \nabla \hat b_v|_y) |  \quad \hat b_v (y) = 0  \}
$$
be the corresponding stable leave. Clearly $\Sigma_v = \pi (\mathcal{H}_v  )$ and since $M$ is of class $C^3$, $\mathcal{H}_v $
is a $C^2$-smooth embedded disc in $S\tm$. As $v$ varies, the sets  $ {\mathcal H}_v $ form a continuous lamination   of $S\tm$ (cf [HI], Proposition 3.1).

Furthermore, since all geodesic balls are convex, the sup-level set $\hat b_v^{-1}([c, \infty))$ is convex for
all $c \in \mathbb R$, see [BGS].

Suppose that  sectional curvatures $K$ of $\tm$ satisfy $ - 1 \le K \le 0$. The standard Hessian comparison
theorem [Pe] asserts
$$
   \| X \|^2 \ge \text{Hess}( -  \hat b_v )(X, X) \ge 0.
$$

\begin{prop}\label{horo} Suppose that  sectional curvatures $K$ of $\tm$ satisfy $ - 1 \le K \le 0$.
Then for any given $\varepsilon >0$, there is $\eta > 0$ such that
if $v' \in {\mathcal H}_v $ satisfies $d_{{\mathcal H}_v} (v,v') < \eta$, then
$ d_{S\tm}(\varphi_t(v), \varphi_t (v')) < \varepsilon $ for all $t \ge 0$.
\end{prop}
\begin{proof} If $\Omega$ is a convex subset of $\tm$,  the nearest point projection
$\mathcal{P}_\Omega: \tm \mapsto \Omega$  is a distance
non-increasing map (see [BGS]). Consider  $\Omega_s = \hat b_v^{-1}([s, \infty))$. Since  $\Omega_s $ is
convex,
$\mathcal{P}_{\Omega_s}$ is distance non-increasing map, and
$$
d_{\Sigma_{\varphi_t(v)}}(\sigma_v(t) , \sigma_{v'}(t)) \le d_{\Sigma_{v}} (v, v') < \eta
$$
for all $t\ge 0 $.

Recall that $\text{Hess}(   \hat b_v )(X, Y) = \langle \nabla_X (\nabla \hat b_v), Y \rangle$. Let
$\Psi_t: [0, \eta] \mapsto \Sigma_{\varphi_t(v)}  $ be a length-minimizing geodesic from
$\sigma_v(t)$ to $\sigma_{v'}(t)$ with respect to the induced metric on the horosphere
$\Sigma_{\varphi_t(v)} = b^{-1}_v(t)$ of height $t$. By the
fact that $\| X \|^2 \ge \text{Hess}( -  \hat b_v )(X, X) \ge 0$ and since
$d_{\Sigma_{\varphi_t(v)}}(\sigma_v(t) , \sigma_{v'}(t))   < \eta $, we obtain, by  integrating along the curve
$\Psi_t $,
$$
     \| \sigma'_{v'}(t) - \mathbb P_{\Psi_t} [ \sigma '_v(t)  ]  \|
     = \| \nabla \hat b_v |_{\sigma_{v'}(t)  }  - \mathbb P_{\Psi_t} [ \nabla \hat b_v |_{\sigma_{v}(t)  }  ]\|
     \le \int_{\Psi_t} \| \nabla ( \nabla \hat b_v  ) \| < \eta,
$$
where $\mathbb P_{\Psi}$ is the parallel translation along the curve $\Psi_t$. It follows that
$$  d_{S\tm}(\varphi_t(v), \varphi_t (v')) < 2 \eta.$$
This completes the proof.
\end{proof}

\

Fix $x\in \tm$ and let $\nu $ be the Patterson-Sullivan  measure
on $\tm (\infty)$ associated to $x$ (see \cite{K1}).
Recall that  $(\tm (\infty ) \times \tm (\infty ))^*$ is the set of distinct pairs of points of the geometric boundary $\tm (\infty)$ and that the action of the covering group $\G$ extends to $\tm (\infty ) $ by continuity and  on $(\tm (\infty ) \times \tm (\infty ))^*$ by $\gamma (\eta, \xi ) = (\gamma \eta, \gamma \xi)$.  In \cite {K1}, \cite {K2}  the following properties of $\nu $ are shown:
\begin {itemize}
\item 1)(\cite {K1}, Lemma 4.1 ) The support of the measure  $\nu\times \nu $ is $(\tm (\infty ) \times \tm (\infty ))^*$.
\item 2)(\cite {K2}, Corollary 4.4) For $(\nu\times \nu ) $ almost every $(\eta, \xi)$, there is a unique regular geodesic $\sigma _{\eta,\xi}$ such that:
$\sigma _{\eta,\xi}(-\infty ) = \eta, \sigma _{\eta,\xi} (+\infty ) = \xi.$
\item 3)(\cite {K2}, Lemma 2.4 ) There is a positive continuous function $F$ on $(\tm (\infty ) \times \tm (\infty ))^*$ such that the measure $\tilde \nu = F   (\nu\times \nu) $ is $\G$ invariant.
\end {itemize}

To a vector $v \in SM$, one associates $Q(v) \in (\tm (\infty ) \times \tm (\infty ))^* \times \mathbb R$ by:
$$ Q(v)\; = \; (\sigma_v (-\infty ), \sigma_v (+\infty ), b(v))$$ where $  b(v) = \lim _{t \to \infty } (d (x, \sigma _v(t)) -t). $
The map $Q$ is a bijection from $S\tm $ on its image. By 2), its image has full $\tilde \nu \times dt $ measure. The measure $ (Q^{-1})_* (\tilde \nu \times dt)$ is therefore a measure on $S\tm$. By 3)  it  is a $\G$ invariant  measure. By definition, it is also   invariant under the  geodesic flow. It corresponds  to a measure $\bar \nu$ on $SM$, which is invariant under the geodesic flow. By 1), the support of $\bar \nu $ is $SM$. By \cite {K2}, Theorem 4.3,   the measure $\bar \nu $ is ergodic under the geodesic flow.

The unit sphere $S_x\tm$ is transversal to the foliation ${\mathcal H}$ and
to  the orbits of the  geodesic flow, so that a tubular neighborhood of $S_x\tm $ will
contain a neighborhood of the  form
$\cup_{v \in S_x\tm} \left(\cup_{s, |s|\le \rho} \varphi _s {\mathcal U}_v\right)$
where ${\mathcal U}_v $ is a neighborhood of $v$ in ${\mathcal H}_v$.
Consider the measure $\nu _x$ on $S_x \tm $ defined by $\nu_x = (P_x^{-1})_* \nu$. On a
neighborhood of $S_x\tm$ of the above form,
the measure $\tilde \nu \times dt$ has a positive density with
respect to the integral over $\nu _x$ of positive measures with
full support on $\cup_{s, |s|\ge \rho} \varphi _s {\mathcal U}_v$
(see \cite {L}, section 3,  for the completely analogous case of negative curvature;
another description of this product structure is in \cite {Gu}). This shows the following:

\begin{prop}\label{measure} Let $A$ be a Borel subset of $S_x\tm$ with $\nu _x(A)>0$. Then, for all $\eta >0$, $(\tilde \nu \times dt) ( \widetilde A_\eta) >0$, where
$$  \widetilde A_\eta \; = \; \cup_{v \in A} \left(\cup_{s, |s|\le \eta } \varphi _s B^{{\mathcal H}_v}(v, \eta)\right)$$
and $B^{{\mathcal H}_v}(v, \eta)$ is the ball of radius $\eta $ in ${\mathcal H}_v$ centered at $v$.
\end{prop}

We now are able to  show that  non-hyperbolic directions are $\nu$ negligible. More precisely, there is

\begin{prop} There exist $h, T $ and $R$ such that, if ${\mathcal F}_K$ is the set of directions $v \in S_x\tm$ such that $\sigma _v (\cdot -t)$ never admits a
$(h,T,R)$ barrier for any $t \ge K$, then  ${\mathcal F}_K$ has no interior in $S_x\tm$ and $\nu_x ({\mathcal F}_K) = 0$.
\end{prop}
\begin{proof}
Recall the set ${\mathcal O}_K$ from section 2. For the sake of the proof, we introduce a slightly smaller set ${\mathcal O}'_K$ which is also generic and full measure, but is disjoint from $(\widetilde {{\mathcal F}_K})_\eta$. The conclusion follows then from Proposition \ref{measure}.
Fix $\delta >0$ small.
\begin{df}\label{bar-delta} We say that the geodesic $\sigma $ admits a $(h,T,R,\delta) $ barrier if there exist $t_i, i= 1, 2,\dots, 6$ with $T_1 +i\delta < t_{i+1} -t_i< T_1 + R - i \delta $  and $t_3 + T <0 < t_4$ such that the geodesic $\sigma $ is $(h,T, {\pi}/{2}+ \delta)$-non flat at  $t_i$, for $i = 1,2,\dots, 6$.
\end{df}
As before, one can find $h,T $ and $R$ such that there is an axis that admits a $(h,T,R,\delta)$ barrier.
 By Proposition \ref{nonflat-cont} the set  ${\mathcal O}'$ of $v$  such that $\sigma _v$ admits  a $(h,T,R,\delta)$ barrier is open. Since the measure $\bar \nu$ is ergodic and has full support, for all positive $K$ the set ${\mathcal O}'_K$ of $v \in SM$ such that the geodesic ray $\sigma _v([K,\infty))$ intersects ${\mathcal O}'$ is open dense in $SM$ and has full $\bar \nu$ measure. By Proposition \ref{nonflat-cont} again, we can find a number $\varepsilon >0$ such that whenever $w \in S\tm$ is such that $\sigma _w$ admits a $(h,T,R,\delta )$ barrier and $d_{S\tm}(w,w') < \varepsilon $, then $\sigma _{w'}$ admits a $(h,T,R)$ barrier. Choose  $\eta$ associated to $\varepsilon $  by Proposition \ref{horo}. Now, if $v \in {\mathcal F}_K$ and $v' \in \cup_{s, |s|\le \eta } \varphi _s B^{{\mathcal H}_v}(v, \eta)$, then $v'$ cannot belong to ${\mathcal O}'_K$ since it would mean that there is a $t>K$ such that $\sigma _{v'}(t) $ admits a $(h,T,R,\delta)$ barrier. Since $d_{S\tm} (\sigma _v(t), \sigma_{v'}(t) ) < \varepsilon$, $\sigma _v(t)$  would admit a $(h,T,R)$ barrier, contrarily to the definition of  ${\mathcal F}_K$. We have shown that  $(\widetilde {{\mathcal F}_K})_\eta $ is disjoint from  ${\mathcal O}'_K$, an open set of full $\bar \nu $ measure. By Proposition \ref{measure},  $\nu_x ({\mathcal F}_K) = 0 $. It is also easy to see that for the same reason,  ${\mathcal F}_K$ has no interior.
\end{proof}

This proves the first part of Theorem \ref{many-mart}, since the set of non-hyperbolic geodesics  starting from $x$ is exactly the union over $K \in {\mathbb N}$ of the  ${\mathcal F}_K$s. For the second part, recall from the introduction the definition of a geodesic ergodic measure on $\tm (\infty)$:
\begin {itemize}
\item 1) The support of the measure  $\mu\times \mu $ is $(\tm (\infty ) \times \tm (\infty ))^*$.
\item 2) For $\mu \times \mu $ almost every $(\eta, \xi)$, there is a unique  geodesic $\sigma _{\eta,\xi}$ such that: $\sigma _{\eta,\xi}(-\infty ) = \eta, \sigma _{\eta,\xi} (+\infty ) = \xi$, and  $\sigma _{\eta,\xi}$ is rank 1.
\item 3) The measure $\mu \times \mu $ is $\G$ quasi-invariant and ergodic: the diagonal action of   $\G$ preserves the $(\mu \times \mu )$negligible subsets  of  $(\tm (\infty ) \times \tm (\infty ))^*$ and  measurable subsets of $(\tm (\infty ) \times \tm (\infty ))^*$ which are $\G$ invariant are either negligible or conegligible.
\end {itemize}

 For $(\eta , \xi  ) \in   (\tm (\infty ) \times \tm (\infty ))^*$, define $N(\eta, \xi) $ as the number of   times, separated by at least $4T_2$,  that the geodesic $\sigma _{\eta, \xi}$,  if it is unique, admits a $(h,T,R)$ barrier.   By property 2)  above the  function $N(\eta , \xi)$ is $(\mu \times \mu) $ almost everywhere well defined. Moreover,  the  function $N(\eta , \xi)$ clearly is $\G$ invariant and therefore  $(\mu \times \mu) $ almost everywhere constant.  We claim that this constant cannot be a finite $K$. Indeed, we just proved that   there is  an open set $\mathcal O''$, such that for $\xi \in \mathcal O''$, there is a regular geodesic $\sigma _\xi $ with $\sigma (0)= x, \sigma (+\infty ) = \xi$ and  at least $K+1$ instants $t_1, \dots, t_{K+1}$ with $t_j -t_{j+1} > 4T_2$,  when  $\sigma _\xi (\cdot - t_j) $ admits  a $(h,T,R)$ barrier . For such a $\xi$, we can find, by (\cite {Ba1})  a small neighborhood  $O_\xi$ of $\sigma _\xi (- \infty)$ such that for $\eta \in O_\xi$, there is a unique $\sigma _{\eta, \xi}$, and it is  close enough to $\sigma _\xi$ that we still have $N(\eta, \xi ) \ge K+1$. Since $\mu \times \mu  \left(\cup _{\xi \in {\mathcal O}''} (O_\xi \times \{\xi\})\right) >0$, this is a contradiction.

So, for $(\mu\times \mu) $ almost every $(\eta, \xi)$, $N(\eta, \xi) $ is infinite.  Let $N_+,N_-,N$ be the subsets of $\{(\eta,\xi): N(\eta, \xi ) = \infty \}$ where there are an infinite number of barrier times respectively only on the positive side  of $\mathbb R$, only on the negative side or on both sides. These three sets are disjoint and $\G$ invariant. Only one of them is of full measure. By Remark \ref {reverse}, the sets $N_+$ and $N_-$ have the same measure, which has to be 0. Therefore, the set $N$ has full measure.    In other words,  $(\mu \times \mu)$ almost every geodesic is hyperbolic. It follows that for  $\mu  $ almost every $\xi \in \tm (\infty)$, there is at least one  geodesic which is asymptotic to $\xi$ and hyperbolic. Theorem \ref {many-mart} follows from Theorem \ref{theo-axis}.


\begin{thebibliography}{An}
\bibitem[An]{An}  Ancona, A., Negatively curved manifolds,
elliptic operators, and the Martin boundary, Ann. Math. vol.
125 (1987), 495-536.

\bibitem[AS]{AS} Anderson, M. and Schoen, R., Positive harmonic functions on complete manifolds of negative curvature, Ann. Math. vol. 121 (1985), 429-461.

\bibitem[Ba1]{Ba1} Ballmann, W., Axial isometries of manifolds
of non-positive curvature, Math. Ann. Vol.259 (1982), 131-144.

\bibitem[Ba2]{Ba2} Ballmann, W., Nonpositively curved manifolds of higher rank. Ann. Math. vol. 122 (1985), 597-609.

\bibitem[Ba3]{Ba3} Ballmann, W., On the Dirichlet problem at infinity for manifolds of nonpositive curvature, Forum Math. vol.1 (1989) 201-213.

\bibitem[Ba4]{Ba4} Ballmann, W., The Martin Boundary of certain Hadamard Manifolds,  Proceedings on Analysis and Geometry (Russian) (Novosibirsk Akademgorodok, 1999),  36--46, Izdat. Ross. Akad. Nauk Sib. Otd. Inst. Mat., Novosibirsk, 2000.

\bibitem[BGS]{BGS} Ballmann, W.,  Gromov, M. and  Schroeder, V.,  Manifolds of nonpositive
curvature. Progress in Mathematics, vol. 61. Birkh\"auser Boston,
Inc., Boston, MA, 1985.

\bibitem[BL]{BL} Ballmann, W. and Ledrappier, F., The Poisson boundary for rank one manifolds and their cocompact lattices, Forum Math. vol.6 (1994) 301-313.

\bibitem[BS]{BS} Burns, K. and Spatzier, R. Manifolds of nonpositive curvature and their buildings. Pub. math. IH\'ES, vol 65, (1987) 35--59


\bibitem[DoC]{DoC} Do Carmo, M., Differential Geometry of Curves
and Surfaces, Prestice-Hall, Inc., New Jersey, 1976


\bibitem[EO]{EO} Eberlein, P. and O'Neill, B., Visibility
manifolds, Pacific J. Math. Vol.46(1973)45-109

\bibitem[GJT]{GJT} Guivarc'h, Y., Ji, L. and Taylor, J.C., Compactifications of symmetric spaces, Birkh\"auser, Boston, 1998

\bibitem[Gu]{Gu} Gunesch, R., Precise asymptotics for the periodic orbits of the geodesic flow in negative curvature, preprint.

\bibitem[HI]{HI} Heitze, E and Im Hof, H.-C., Geometry of horospheres, J. Diff. Geom. vol.12 (1977) 481-491.

\bibitem[K1]{K1} Knieper, G., On the asymptotic geometry of nonpositively curved manifolds, GAFA, vol.7 (1997), 755--782.

\bibitem[K2]{K2} Knieper, G., The uniqueness of the measure of maximal entropy on rank 1 manifolds, Ann. math. vol.148 (1998), 291-314.

\bibitem[Ka]{Ka} Kaimanovich, V. A.,  Boundaries of invariant Markov operators: the identification problem,
Ergodic theory of $Z^ d$ actions (Warwick, 1993--1994), 127--176,
London Math. Soc. Lecture Note Ser., 228,


\bibitem[L]{L} Ledrappier, F., Harmonic measures and Bowen-Margulis measures, Israel J. Math. vol.71 (1990) 275-282.

\bibitem[MV]{MV} Mazzeo, R. and Vasy, A., Resolvents and Martin boundary of product spaces, GAFA, vol.12 (2002) 1018--1079


\bibitem[Pe]{Pe} Petersen, P., Riemannian Geometry,
Springer-Verlag, GTM. Vol.171, New York 1998

\bibitem[SY]{SY} Schoen, Y. and Yau, S.-T., Lectures on Differential Geometry, International Press, Boston, 1994.



\end{thebibliography}
\end {document}